\documentclass[10pt,twoside]{article}
\usepackage{rotating,graphicx,amssymb,amsmath}
\usepackage{pst-all}

\newcommand{\bs}{\boldsymbol}
\newcommand{\pd}[2]{\dfrac{\partial #1}{\partial #2}}
\newcommand{\vlegend}[3]{\turnbox{90}{\mbox{ }\hspace{#1 in}\small #2} \hspace{-#3in} }

\begin{document}

\pagestyle{myheadings}\markboth{T. Dubois, R. Touzani}{Heat island flows}

\title{A numerical study of heat island flows in an open domain: Stationary solutions}

\author{Thierry Dubois\thanks{Email: Rachid.Touzani@univ-bpclermont.fr} and
        Rachid Touzani\thanks{Email: Thierry.Dubois@math.univ-bpclermont.fr} \\
    Laboratoire de Math\'ematiques, UMR CNRS 6620, \\
    Universit\'e Blaise Pascal (Clermont--Ferrand 2), \\
    63177 Aubi\`ere cedex, France}

\maketitle

\begin{abstract}
We present two dimensional numerical simulations of a natural convection problem in an unbounded domain. 
A thermal stratification is applied in the vertical direction and the flow circulation is induced by
a heat island located on the ground.
For this problem, thermal perturbations are convected in the horizontal direction far from the heated element
so that very elongated computational domains have to be used in order to compute accurate numerical solutions. 
To avoid this difficulty thermal sponge layers are added at the vertical boundaries.
With this approach, stationary solutions at $Ra\le 10^5$ are investigated. 
Boussinesq equations are discretized with a second-order finite volume scheme on a staggered grid combined with 
a second-order projection method for the time integration. 
\end{abstract}


\section{Introduction}
In this paper we present numerical simulations of a particular type of thermal fluid flows.
Namely, we are concerned with the so-called \textit{heat island} flows \cite{MalkusStern}
\textit{i.e.} fluid flows where natural convection is generated by a local variation of temperature thus inducing
buoyancy effect. This phenomenon appears in the presence of heat stratification
that stabilizes the fluid flow. The present model is generally used to study environment problems
such as urban heat island \cite{OlfeLee,Delage-Taylor}. 
Heat island fluid flows occur in open configurations, which require, from the mathematical
viewpoint, their study in unbounded domains. In practice, numerical simulations are carried out
in \textit{large} but bounded computational domains for which appropriate design of boundary conditions 
must be investigated. This issue is a central one in our contribution.

For this problem, the heat island perturbation generates an ascending flow circulation
which develops mainly in an area surrounding the heated element. The vertical stratification limits this effect 
by pushing the flow down to the ground. 
As a consequence of these opposite forces, thermal perturbations are propagated in the horizontal direction
at long distance far from the heat source. Therefore, \textit{very elongated} domains have to be used 
in order to accurately compute the temperature deviation from the stratified profile.

Despite the increase of computational resources, the direct numerical simulation of solutions to natural convection 
induced by a local heat source in large domains remains a real challenge (see \cite{Xin-al} and the
references therein). In most cases, the far-field solutions are unknown so that the use of a limited computational 
area surrounding the heating element requires an appropriate treatment of the boundary conditions. In \cite{Xin-al},
outer artificial conditions are applied at the domain boundaries. 
For heat island flows, only the length of the computational domain has to be limited. Flow elevation being 
strongly reduced by the vertical stratification, the domain length necessary to obtain converged solutions is not
large. The approach used in this paper consists in applying a 
thermal sponge layer in the vicinity of the vertical boundaries: The temperature equation is modified so 
that the convective terms are \textit{smoothly} damped in an area closed to the outflow boundaries.
Boundary conditions with a sponge are classically employed in computational electromagnetics
\cite{Berenger,Abarbanel-al} and in simulations of compressible turbulent flows \cite{Mahesh-al,Boersma}.
We show that this technique is well suited for the numerical simulations of heat island type flows 
at moderate Rayleigh numbers, that is $Ra\le 10^5$. This issue constitutes the main contribution of this paper.

The outline of the paper is as follows. In the next section, we describe the set of equations that govern the 
fluid flow in a heat island as well as the domain geometry and boundary conditions. In particular, we write 
the equations in a nondimensional form that involves two parameters, the Rayleigh number and a thermal 
stratification coefficient. We show that a specific choice of relevant parameters is to be made here for 
the heat island flow. Section 3 describes approximations in closed domains. Two approaches are used: The 
former relies on the use of \textit{very large} computational domains without any particular treatment at
the domain exits while the latter introduces thermal sponge layers acting in the vicinity of the vertical 
boundaries. Section 4 presents the space and time discretization schemes.  
Preliminary numerical results were obtained by Touzani in \cite{Touzani} using a finite element method
coupled with a penalty method to impose the incompressibility constraint.
Here, we use a second-order finite volume scheme on a staggered grid for space 
discretization and a second-order projection method for the time integration of the resulting system of 
differential equations. 
Section 5 gives numerical results. Numerical simulations in a square differentially heated cavity are 
first performed to check the code accuracy. Stationary solutions to the heat island problem at $Ra\le 10^5$ 
are obtained in \textit{very elongated} computational domains and are used as references to validate 
the sponge technique applied to the heat equation. Accurate stationary solutions are then computed with 
this approach: Characteristic values are listed for reference. Finally, conclusions are drawn and perspectives 
for future works on this problem are discussed.

\section{Description of the problem}

\subsection{The physical problem}
We consider a fluid that fills the half plane $\left\{\bs{x}^\star=(x^\star,y^\star)\in\mathbb{R}^2; y^\star>0\right\}$.
Here and in the sequel we shall append a superscript ${}^\star$ to all physical dependent and independent variables, 
the notation without ${}^\star$ being reserved to nondimensional variables. 
The fluid is initially at rest and is thermally stratified in the vertical direction, namely the velocity field
$\bs{u}^\star=(u^\star,v^\star)$ and the potential temperature $T^\star$ satisfy, at time $t^\star=0$,
\begin{align}
  & \bs{u}^\star=\bs{0},     \\[0.15cm]
  & T^\star\,=\,T_0\,+\, \alpha_\textrm{s}\, y^\star, 
  \label{BC-1}
\end{align}
where $T_0>0$ is the potential temperature at the ground and $\alpha_\textrm{s}>0$ is the thermal stratification 
coefficient.

In order to generate a flow, a local temperature perturbation of intensity $T_1>0$ is applied on a source line 
$Q^\star=(-\delta/2,\delta/2),\ \delta>0,$ located on the ground (see Figure 1), that is we impose
\begin{equation} 
  T^*(\bs x^*,t^*) = T_0 + \frac{T_1}{2} 
  \biggl(1-\tanh{\Bigl(\frac{2\vert x^*\vert+\delta}{2\,\delta\,\zeta}\Bigr)}\biggr),   \label{BC-ground^*} 
\end{equation}
for all $\bs{x}^*\in \mathbb{R}^2\cap\{y^*=0\}$ and for all time $t^*>0$. Note that this thermal perturbation is 
constant in time and is a regularized version of a heat island type perturbation (see \cite{Delage-Taylor} for instance) 
for which a constant and uniform temperature would be applied on the heated element $Q^\star$.
The parameter $\zeta>0$ in \eqref{BC-ground^*} is used to set the sharpness of the temperature gradient 
$\partial T^*/\partial x^*$ near the plate boundaries $\vert x^*\vert=\delta/2$ at the ground level $y^*=0$.
In this study, the value $\zeta=2.5\times 10^{-2}$ is used.

Due to the perturbation \eqref{BC-ground^*}, a thermal plume develops above the heated plate $Q^\star$.
Natural convection induces an ascending flow circulation while the gravity force and the vertical stratification 
limit the development of flow structures in the vertical direction.
As a consequence of these opposite forces, thermal perturbations are propagated in the horizontal direction
at long distance, far from the heated element. 
We may then expect that the solutions will decay \textit{rapidly} in the vertical direction and \textit{very slowly} 
in the horizontal one.
The main difficulty for this problem, as long as numerical simulations are concerned, resides in a suitable design
of boundary conditions in order to properly reproduce the behavior of the far-field solutions.
Indeed, errors in the numerical approximation at long distance from the heated plate may deteriorate 
the accuracy of simulations in the region surrounding the heated element.

\subsection{The governing equations}

We consider the set of equations describing a two-dimensional thermal flow assuming Boussinesq approximation.
Let us recall that this one stipulates that for small temperature differences, the density variations
are more significant in the gravity acceleration term than in others. 

Velocity $\bs{u}^*$, pressure $p^*$, density $\rho^*$ and potential temperature $T^*$ satisfy
the set of equations:
\begin{align}
   & \pd{\bs{u}^*}{t^*} - \nu\,\Delta^* \bs{u}^* + \bs{\nabla}^*\cdot(\bs{u}^*\otimes\bs{u}^*) 
                    + \frac 1{\rho_0}\bs{\nabla}^* p^* = -\frac{\rho^*}{\rho_0}g\,\bs{e}_2,     \label{Boussinesq^*-1} \\
   & \bs{\nabla}^*\cdot\bs{u}^* = 0,                                                            \label{Boussinesq^*-2} \\
   & \pd{T^*}{t^*} - \kappa\,\Delta^*T^* + \bs{\nabla}^*\cdot(\bs{u}^*T^*) = 0,                 \label{Boussinesq^*-3} \\
   & \bs{u}^*(\bs x^*,0) = \bs{0},\quad T^*(\bs{x}^*,0) = T_0+\alpha_s\, y^*,                   \label{Boussinesq^*-4}
\end{align}
where the physical constants are the kinematic viscosity $\nu$, the thermal conductivity $\kappa$
and the modulus of gravity acceleration $g$. The unit vector in the vertical direction is denoted
by $\bs{e}_2$, namely $\bs{e}_2=(0,1)$.

For these type of flows, we have adopted the approximation that compressibility is expressed by
a dilatation equation. Therefore, density is related to the temperature variations with respect
to the reference state $(T_0,\rho_0)$, by the following equation
\begin{equation}
  \rho^* = \rho_0(1-\beta(T^*-T_0)),   \label{Boussinesq^*-5}
\end{equation}
where $\beta$ is the thermal expansion coefficient.
Equations \eqref{Boussinesq^*-1}--\eqref{Boussinesq^*-5} are to be considered in the infinite domain 
$\mathbb R\times\mathbb R^+$ where the line $y^*=0$ contains the heated plate $Q^*=(-\delta/2,\delta/2)$. 
Such a problem with appropriate behavior at the infinity is well suited for generating a heat island flow in the 
vicinity of the heated plate if a thermal stratification is given.

The boundary condition \eqref{BC-ground^*} has as effect to generate, in the neighborhood of the plate $Q^*$ a
temperature plume with a shape and an intensity depending on the system parameters. 

\subsection{Nondimensional form of the equations}

In order to normalize the problem, we introduce as usual reference values for temperature
$T_r=T_1$, length $L_r=\delta$ and velocity $U_r=\sqrt{g\beta L_r T_r}$. 
Reference values for time and pressure can then be deduced by $t_r=L_r/U_r$ and $p_r=\rho_0\,U_r^2$ 
respectively. We define nondimensional space and time variables by $\bs x=\bs x^*/L_r$ and $t=t^*/t_r$.
In terms of these nondimensional variables, the heated plate reads $Q=(-1/2,1/2)$. 
The nondimensional variables
\begin{equation}
   \bs{u}=\frac{\bs{u}^*}{U_r}, \ T=\frac{T^*-T_0}{T_r} \quad \textrm{and}\quad p=\frac{p^*+\rho_0\, g\,y^*}{p_r}
   \label{nondim-var}
\end{equation}
satisfy in $\mathbb{R}\times\mathbb{R}^+$ and for $t>0$ the following system of equations:
\begin{align}
   & \pd{\bs u}{t}-\sqrt{\frac{Pr}{Ra}}\,\Delta\bs u + \bs{\nabla}\cdot(\bs u\otimes\bs u) 
                                                   + \bs{\nabla} p = T\,\bs{e_2},         \label{Boussinesq-1} \\[0.2cm]
   & \bs{\nabla}\cdot\bs{u} = 0,                                                          \label{Boussinesq-2} \\[0.2cm]
   & \pd{T}{t} - \frac{1}{\sqrt{Ra\,Pr}}\,\Delta T + \bs{\nabla}\cdot(\bs{u}\,T) = 0,   \label{Boussinesq-3} \\[0.2cm]
   & \bs{u}(\bs{x},0) = \bs{0},\quad T(\bs{x},0) = \alpha\, y,                            \label{Boussinesq-4}
\end{align}
where $\alpha=\alpha_s\,L_r/T_r$. The Prandtl and Rayleigh numbers are respectively defined by
$$ Pr = \frac\nu\kappa\quad\textrm{and}\quad Ra=\frac{g\,\beta\,L_R^3\,T_R}{\nu\,\kappa}. $$
The nondimensional form of the boundary condition \eqref{BC-ground^*} reads
\begin{equation}
  T(\bs{x},t) = \frac{1}{2} \biggl(1-\tanh{\Bigl(\frac{2\vert x\vert+1}{2\,\zeta}\Bigr)}\biggr)\quad\textrm{for}\ \ 
                          \bs{x}=(x,0). 
  \label{BC-ground} 
\end{equation}
To simplify, we set the Prandtl number, that characterizes the fluid to its value for the air,
$Pr=0.71$. Therefore, the system of equations \eqref{Boussinesq-1}--\eqref{Boussinesq-4} depends on
two parameters $Ra$ and $\alpha$. 

The aim of this work is to study through \textit{accurate} numerical simulations the behavior of the stationary
solutions and their dependency on the Rayleigh number, at a fixed stratification coefficient $\alpha$.

\section{Approximation in closed domain}

\subsection{Large elongated domain}

Clearly, Problem \eqref{Boussinesq-1}--\eqref{BC-ground} is difficult to handle numerically in an unbounded domain.
We choose to approximate $\mathbb{R}\times\mathbb{R}^+$ by a rectangle $\Omega=(-L/2,L/2)\times(0,H)$ with $L$ large
enough (see Figure \ref{fig:comp-domain}) and we denote by $\Gamma$ its boundary.
A simple and naive approach would consist in imposing conditions on $\Gamma$ consistent with the initial condition
\eqref{Boussinesq-4} and the boundary condition \eqref{BC-ground}. Therefore, equations 
\eqref{Boussinesq-1}--\eqref{Boussinesq-4} are supplemented with the boundary conditions:
\begin{alignat}{2}
   & \bs u(\bs x,t) = 0,              && \qquad\bs x\in\Gamma,                               \label{bc-1} \\
   & T(\bs x,t) = \frac{1}{2} \biggl(1-\tanh{\Bigl(\frac{2\vert x\vert+1}{2\,\zeta}\Bigr)}\biggr),
                                      && \qquad\bs x\in\Gamma_0=\{\,\bs x\in\Gamma;\ y=0\},  \label{bc-2} \\
   & T(\bs x,t) = \alpha\, y,         && \qquad\bs x\in\Gamma\setminus\Gamma_0,              \label{bc-3}
\end{alignat}
for all time $t>0$. 

Such conditions can be enforced only if the domain lengths $L$ and $H$ are large enough
so that their change has a negligible effect on the solution.
As it was previously mentioned, the gravity force and the stratified profile \eqref{bc-3} tend to push the flow 
down to the ground, limiting the vertical convection. 
Therefore, the flow variables have a rapid decay with respect to the elevation $y$, so that the domain height $H$ does not 
need to be too large. 
The heat island perturbation generates an ascending flow circulation which is essentially local. On the other hand, 
thermal perturbations are convected in the horizontal direction at long distance 
far from the heated source line. Therefore, very elongated computational domains have to be considered, 
that is $L\gg H$ and $L\gg 1$, in order to produce accurate solutions. 

\subsection{A truncated temperature equation}

We shall see in the development of this study that the above approach is very consuming in terms of computational 
resources inducing strong limitations on the allowable grid resolution. Indeed, if the domain length $L$ is not
large enough, artificial boundary layers develop at the domain boundaries $\vert x\vert=L/2$ inducing an 
overestimation of the flow variables. These numerical errors deteriorate the accuracy of the solutions in the 
central area where most of the flow dynamics take place.  
An appropriate design of the behavior of the solution close to these boundaries is necessary in order to relax
the condition $L\gg 1$ on the domain length.
To do so, we propose to limit the horizontal propagation of the thermal perturbation by damping
the convective terms in the temperature equation in sponge layers close to the domain exits $\vert x\vert=L/2$. 
The nonlinear convection term in \eqref{Boussinesq-3} is multiplied by a filter function 
\begin{equation} \label{psi-gamma}
  \psi_\gamma(x)=e^{-\gamma\left(\frac{2|x|}{\sigma L}\right)^p} 
\end{equation}
where $\sigma\in(0,1)$ and $p\ge 1$. This yields the modified heat equation
\begin{equation}
  \pd{T}{t} - \frac{1}{\sqrt{Ra\,Pr}}\,\Delta T + \psi_\gamma(x)\,\bs{\nabla}\cdot(\bs{u}\,T) = 0.
  \label{Boussinesq-3-filt} 
\end{equation}
For the sake of simplicity, we use the same notation $T$ for the \textit{truncated} and \textit{standard} temperature 
respectively solution of \eqref{Boussinesq-3-filt} and \eqref{Boussinesq-3}. 
The former corresponds to the choice of $\gamma=1$ and the latter to $\gamma=0$.

The filter function $\psi_\gamma$ rapidly decays when $|x|\approx L/2$ whereas $\psi_\gamma\longrightarrow 1$ in the
center of the computational domain thus reducing to a classical convective term. 
The effect of $\psi_\gamma$ is to introduce a \textit{sponge layer}, close to the vertical boundaries, where the 
convection of temperature is smoothly damped through the outflow.
As we will see in Section 5, such a treatment allows to significantly reduce the size of the computational domain 
$\Omega$ required to reach a given accuracy.

This approach, while different in its implementation, is similar to the perfectly matched layer (PML) method used in 
computational electromagnetics and introduced by Berenger \cite{Berenger}. Boundary conditions with a sponge are also
commonly used for the numerical simulations of compressible turbulent flows as jet flows for instance \cite{Boersma}.
For these problems, the solution is driven to a specified outflow state by adding in the sponge layer a
\textit{cooling} term to the right-hand side of the equations. 
We found that our method, for the heat island problem \eqref{Boussinesq-1}--\eqref{BC-ground}, is less sensitive to
the parameters involved in the definition of the sponge functions.

\subsection{A model formulated in terms of temperature fluctuations}
By noting that a stratified profile of the temperature can be expressed in the momentum equation as a gradient term,
we decompose the potential temperature $T$ into $T=\alpha\, y+\theta$ which introduces the temperature fluctuation
$\theta$. By reporting this decomposition into \eqref{Boussinesq-1}, \eqref{Boussinesq-2}, \eqref{Boussinesq-3-filt}
and \eqref{Boussinesq-4}, recalling that $\bs u=(u,v)$ and, introducing the new pressure variable $P=p-\alpha\,y^2/2$, 
we finally obtain the system of equations:
\begin{align}
   & \pd{\bs u}{t}-\sqrt{\frac{Pr}{Ra}}\,\Delta\bs u + \bs{\nabla}\cdot(\bs u\otimes\bs u) 
                                                   + \bs{\nabla} P = \theta\,\bs{e_2},    \label{Boussinesq-5} \\[0.2cm]
   & \bs{\nabla}\cdot\bs{u} = 0,                                                          \label{Boussinesq-6} \\[0.2cm]
   & \pd{\theta}{t} - \frac{1}{\sqrt{Ra\,Pr}}\,\Delta \theta 
           \,+\, \psi_\gamma(x)\,\bigl(\bs{\nabla}\cdot(\bs{u}\,\theta)\,+\,\alpha\,v\bigr)= 0,  \label{Boussinesq-7} \\[0.2cm]
   & \bs{u}(\bs{x},0) = \bs{0},\quad \theta(\bs{x},0) = 0,                                \label{Boussinesq-8}
\end{align}
which is supplemented with:
\begin{alignat}{2}
   & \bs u(\bs x,t) = 0,              && \qquad\bs x\in\Gamma,                               \label{bc-4} \\
   & \theta(\bs x,t) = \frac{1}{2} \biggl(1-\tanh{\Bigl(\frac{2\vert x\vert+1}{2\,\zeta}\Bigr)}\biggr),
                                      && \qquad\bs x\in\Gamma_0=\{\,\bs x\in\Gamma;\ y=0\},  \label{bc-5} \\
   & \theta(\bs x,t) = 0,                   && \qquad\bs x\in\Gamma\setminus\Gamma_0.        \label{bc-6}
\end{alignat}
In the next section, the numerical approximation of \eqref{Boussinesq-5}-\eqref{bc-6} is addressed.

\section{Numerical approximation}

The numerical discretization of \eqref{Boussinesq-5}--\eqref{bc-6} is achieved by using a second-order
projection scheme in time coupled with a second-order finite volume approximation in space. The unknowns are placed 
on a staggered mesh as for the classical MAC scheme \cite{HarlowWelch}.

\subsection{Time discretization}

The natural convection problem \eqref{Boussinesq-5}--\eqref{bc-6} is solved in two steps decoupling the computation 
of the temperature fluctuation and of the velocity-pressure unknowns. A second-order projection scheme 
\cite{Gresho,Guermond,Quarteroni} is first applied to solve the momentum equations
\eqref{Boussinesq-5} and to enforce the incompressibility constraint \eqref{Boussinesq-6}.

Let $\delta t>0$ stand for the time step and $t^k=k\,\delta t$ discrete time values. 
Let us consider that $(\bs u^j,\, P^j,\,\theta^j)$ are known for $j\le k$. The computation of 
$(\bs u^{k+1},\, P^{k+1})$ consists in:

-- Computing a predictor $\widetilde{\bs u}^{k+1}$ by solving:
\begin{align}
   & \dfrac{\widetilde{\bs u}^{k+1}-{\bs u}^k}{\delta t} -
            \sqrt{\frac{Pr}{Ra}}\,\Delta\Big(\dfrac{\widetilde{\bs u}^{k+1}+\bs u^k}2\Big)
                                          + \bs\nabla P^k = \frac 12\,(3\,\theta^k-\theta^{k-1})\,\bs e_2  \nonumber\\[0.2cm]
   & \qquad\qquad\qquad\qquad - \frac 32\,\bs\nabla\cdot(\bs u^k\otimes\bs u^k) 
                                     + \frac 12\,\bs\nabla\cdot(\bs u^{k-1}\otimes\bs u^{k-1}),   \label{Boussinesq-1a-dt}\\[0.2cm]
   & \widetilde{\bs u}^{k+1} = \bs 0\quad\text{on }\Gamma.                                          \label{Boussinesq-2a-dt}
\end{align}

-- Projecting to obtain a divergence free velocity $\bs u^{k+1}$:
\begin{align}
   &\dfrac{\bs u^{k+1}-\widetilde{\bs u}^{k+1}}{\delta t} + \frac 12\,\bs\nabla(P^{k+1}-P^k) = 0,  \label{Boussinesq-1b-dt}\\[0.2cm]
   &\bs\nabla\cdot\bs u^{k+1} = 0, \quad \bs u^{k+1}\cdot\bs n = 0 \quad \text{on }\Gamma.       \label{Boussinesq-2b-dt}
\end{align}
Finally, the temperature variation $\theta^{k+1}$ is computed by solving:
\begin{align}
   & \dfrac{\theta^{k+1}-\theta^{k}}{\delta t} 
           - \frac 1{\sqrt{Ra\,Pr}}\,\Delta\Big(\frac{\theta^{k+1}+\theta^k}2\Big)              \nonumber\\[0.2cm]
   & \quad = -\psi_\gamma(x) \left(\frac 32\,\bigl(\bs\nabla\cdot(\bs u^k\theta^k)+\alpha\,v^k\bigr)
                    -\frac 12\,\bigl(\bs\nabla\cdot(\bs u^{k-1}\theta^{k-1})+\alpha\,v^{k-1}\bigr)\right), \label{Boussinesq-1c-dt}\\[0.2cm]
   & \theta^{k+1}=\frac 12\,\biggl(1-\tanh{\Bigl(\frac{2\vert x\vert+1}{2\,\zeta}\Bigr)}\biggr)\quad\text{on }\Gamma_0,\quad
     \theta^{k+1}= 0\quad\text{on }\Gamma\setminus\Gamma_0.                                    \label{Boussinesq-2c-dt}
\end{align}

Hence, viscous and diffusion terms are discretized with a Crank-Nicholson scheme while 
nonlinear convective terms are integrated by an Adams-Bashforth scheme. The scheme 
\eqref{Boussinesq-1a-dt}--\eqref{Boussinesq-2c-dt} is globally second-order accurate.
This method is well-suited for the Navier-Stokes equations and is frequently used (see for 
instance \cite{KimMoin} and \cite{LeQuere-DeRoquefort}).

Finally, note that the temperature is computed once the projected (divergence free) velocity is obtained.  
A different approach is used in \cite{LeQuere,Bruneau-Saad,Xin-LeQuere}: the temperature is 
first computed, then the velocity field is obtained by using the temperature at the new time level.
As we will see in Section 5.1, both time stepping schemes achieve a second-order accuracy and provide similar
results.
Also, in \cite{LeQuere} and \cite{Bruneau-Saad}, a BDF projection scheme is implemented to discretize the advection-diffusion
terms. The BDF scheme may be more efficient for the computation of nonstationary solutions as it provides a more
accurate approximation of the pressure. In the case of stationary solutions, which are our main concern in this paper,
the more classical approach used here provides satisfactory results.

\subsection{Space discretization}

\subsubsection{Mesh and unknown locations}

Due to the combined effects of the gravity force and the vertical stratification, flow variables decay rapidly
with respect to the vertical elevation. Therefore, the domain $\Omega$ is discretized by using 
a uniform subdivision in the $y$-direction. A non-uniform grid is required in the $x$-direction as $L\gg H$.
Let $N$ and $M$ denote two integers and let
\begin{alignat*}{3} 
   & x_i = \frac L2\,\varphi(i\,\ell) && \qquad\text{for } i=0,\ldots,N,  && \quad \ell=\frac LN,    \\[0.2cm]
   & y_j = j\,h                     && \qquad\text{for } j=0,\ldots,M,  && \quad    h=\frac HM.
\end{alignat*}
The function $\varphi:(0,L)\longrightarrow (-1,1)$, describing mesh density, is defined by
\begin{equation} \label{map-mesh}
  \varphi(x) = \frac{2x-L+\gamma_1\,\tanh(\gamma_2x)-\gamma_1\,\tanh(\gamma_2(L-x))}{L+\gamma_1\,\tanh(\gamma_2L)}. 
\end{equation}
The function $\varphi$ enables defining a grid with steps $\ell_i=x_i-x_{i-1}$ and
distortion ratios $r_i=\ell_i/\ell_{i-1}$ that increase in function of the distance from the center of
the heated element. The parameters $\gamma_1$ and $\gamma_2$ are chosen so that the lengths $\ell_i$ 
are of order $h$ in the neighborhood of the heated plate $Q=(-1/2,1/2)$.

We introduce points $x_{i+1/2}:=\frac{x_i+x_{i+1}}{2}$ for $i=0,\ldots,N-1$, and
$y_{j+1/2}:=\frac{y_j+y_{j+1}}{2}$ for $j=0,\ldots,M-1$.
All terms in equations \eqref{Boussinesq-1a-dt}, \eqref{Boussinesq-1b-dt}, \eqref{Boussinesq-2b-dt} and 
\eqref{Boussinesq-1c-dt} are discretized in space by using second-order centered finite volume schemes. The discrete 
unknowns are given on a staggered grid (see \cite{HarlowWelch}): discrete pressure values are located at the center 
of mesh cells $K_{i-\frac 12,j-\frac 12}=(x_{i-1},x_i)\times (y_{j-1},y_j),$ vertical velocity and temperature values 
are located at the center of mesh cells $K_{i-\frac 12,j}=(x_{i-1},x_i)\times (y_{j-\frac 12},y_{j+\frac 12}),$
and those of the horizontal velocity are located at the center of mesh cells 
$K_{i,j-\frac 12}=(x_{i-\frac 12},x_{i+\frac 12})\times (y_{j-1},y_j),$ as it is shown on Figure \ref{fig-cell}.
We define the vector $\bs u^{k+1}\in\mathbb{R}^{(N-1)M}$ of components $u^{k+1}_{ij}$
and similarly, $\bs v^{k+1}\in\mathbb{R}^{N(M-1)}, \bs P^{k+1}\in\mathbb{R}^{(N-1)(M-1)}$ and 
$\bs\theta^{k+1}\in\mathbb{R}^{N(M-1)}$ of components $v^{k+1}_{ij}, P^{k+1}_{ij}$ and $\theta^{k+1}_{ij}$ respectively. 

In \cite{Bruneau-Saad}, the discrete temperature is located at the pressure nodes $(x_{i-\frac 12},y_{j-\frac 12})$. 
However, this choice implies interpolations in order to compute the contribution of temperature in the vertical
velocity momentum equation and \textit{vice versa}. A more convenient choice for natural convection problems is
to place temperature at the same nodes as the vertical velocity.

\subsubsection{Discrete systems}

The discretization of \eqref{Boussinesq-1a-dt} is achieved by integration of the equation of the horizontal 
(resp. vertical) velocity component over the volume cells $K_{i,j-\frac 12}$ (resp. $K_{i-\frac 12,j}$). 
Gradient and Laplace operators are classically approximated by centered second-order finite volume schemes. 
Approximation of the nonlinear terms requires second-order interpolations of velocity components at the interfaces, 
for instance we use
$$ \int_{y_{j-1}}^{y_j} u^2(x_{i-\frac 12},y)\,dy  \approx h\,\Bigl(\frac{u^2_{ij}+u^2_{i-1,j}}2\Bigr) $$
and
$$ \int_{x_{i-\frac 12}}^{x_{i+\frac 12}} (uv)(x,y_j)\,dx
             \approx \Bigl(\frac{u_{ij}+u_{i+1,j}}2\Bigr)\, \frac{(\ell_{i+1}v_{ij}+\ell_i v_{i+1,j})}{(\ell_{i+1}+\ell_i)}. $$
Similar interpolation rules are also applied to discretize the equation satisfied by the vertical velocity $v$.
This leads to the system of equations
\begin{align}
   \widetilde{\bs u}^{k+1} + \frac{\delta t}2\,\sqrt{\frac{Pr}{Ra}}\, A_1\widetilde{\bs u}^{k+1}  \nonumber
             & = -\,\delta t\, G_1 \bs P^k + \bs u^k - \frac{\delta t}2\, \sqrt{\frac{Pr}{Ra}}\, A_1 \bs u^k   \\
             &  - \frac{\delta t}2\, \left(3\,N_1(\bs u^k,\bs v^k)-N_1(\bs u^{k-1},\bs v^{k-1})\right),  \label{Num_sch-1}  \\
   \widetilde{\bs v}^{k+1} +\frac{\delta t}2\sqrt{\frac{Pr}{Ra}}\,  A_2 \widetilde{\bs v}^{k+1}
             & = -\, \delta t\,G_2 \bs P^k + \bs v^k - \frac{\delta t}2\,\sqrt{\frac{Pr}{Ra}}\,A_2\bs v^k   \nonumber \\
             &-\frac{\delta t}{2}\left(3\,N_2(\bs u^k,\bs v^k)-N_2(\bs u^{k-1},\bs v^{k-1})\right)  \label{Num_sch-2} \\
             &  + \frac{\delta t}2 (3\,\bs\theta^{k}-\bs\theta^{k-1}),                             \nonumber 
\end{align}

where the matrices $A_i$ are discrete approximations of the operator $-\Delta$ with appropriate treatment of the
boundary conditions for the velocity components, $G_i$ are the ones of the gradient components 
and $N_i$ are the ones of the nonlinear terms.

With the use of the staggered implementation of the discrete values on the mesh, several possibilities are offered
for the treatment of boundary conditions. Concerning the vertical velocity component, we choose to impose the boundary 
conditions on vertical boundaries at grid points $\{(x_0,y_j),j=1,\ldots,M-1\}$ and $\{(x_N,y_j),j=1,\ldots,M-1\}$. 
This yields a modified formula for the discretization of $\frac{\partial^2 v}{\partial x^2}$ at the first point
away from the vertical boundary, that is 
$$\int_{K_{1/2,j}}\frac{\partial^2 v}{\partial x^2}(x,y)\,dx\,dy \approx
     h \,\Biggl(\frac{v_{3/2,j}-v_{1/2,j}}{x_{3/2}-x_{1/2}}-\frac{v_{1/2,j}-v_{0,j}}{x_{1/2}-x_0}\Biggr).$$
A similar formula applies at the last inner point in the horizontal direction, that is $x_{N-1/2}$.
On horizontal boundaries, boundary conditions for $v$ are imposed at points $\{(x_{i-\frac 12},y_0),i=1,\ldots,N\}$ and 
$\{(x_{i-\frac 12},y_M),i=1,\ldots,N\}$.
In the vertical direction, the use of second-order centered formula and of uniform mesh points allows us to
apply a discrete Fourier transform \cite{Schumann-Sweet}. We thus obtain a set of independent and symmetric tridiagonal
systems which can be efficiently solved with the LDL$^\textrm{T}$ algorithm.

Boundary conditions for $u$ on the vertical boundaries are also imposed at mesh points, that is 
$\{(x_0,y_{j-\frac 12}), j=1,\ldots,M\}$ and $\{(x_N,y_{j-\frac 12}), j=1,\ldots,M\}$.
However, at the top and bottom horizontal boundaries, values at \textit{ghost points} $y_{-1/2}=-\frac h2$ and
$y_{M+1/2}=H+\frac h2$ are used to impose boundary conditions with a second-order extrapolation formula: 
We introduce \textit{ghost} velocity values
\begin{equation*} 
  u_{i,-1}=-u_{i,1} \quad\textrm{and}\quad u_{i,M+1}=-u_{i,M} \quad\textrm{for} \quad i=1,\ldots,N.
\end{equation*}
The discretization of $\partial^2/\partial y^2$ on the sequence of mesh points $y_{1/2},\ldots,y_{M-1/2}$ with
a second-order centered finite volume scheme also yields a discrete operator which can be easily diagonalized by applying
a discrete Fourier transform \cite{Schumann-Sweet}.

The discrete version of \eqref{Boussinesq-1b-dt} is obtained similarly: 
\begin{equation} \label{Num_sch-3}
  \begin{split}
    &\bs u^{k+1}=\widetilde{\bs u}^{k+1}-\frac{\delta t}2\,G_1\bs\phi^{k+1}, \\[0.25cm]
    &\bs v^{k+1}=\widetilde{\bs v}^{k+1}-\frac{\delta t}2\,G_2\bs\phi^{k+1},
  \end{split}
\end{equation}
where $\bs\phi^{k+1}=\bs P^{k+1}-\bs P^k.$ 
Note that due to the staggered locations of the unknowns, no boundary conditions for the pressure are required in the 
correction step \eqref{Num_sch-3}. Therefore, discrete pressure is defined only at interior points.

The discretization of the incompressibility constraint is achieved by integrating \eqref{Boussinesq-2b-dt} over the
pressure cell $K_{i-\frac 12,j-\frac 12}$, leading to
\begin{equation} \label{Num_sch-4}
     D_1\bs u^{k+1} + D_2\bs v^{k+1} = 0,
\end{equation}
where $D_1$ and $D_2$ are approximations of $\partial/\partial x$ and $\partial/\partial y$.
Combining \eqref{Num_sch-3} and \eqref{Num_sch-4}, we deduce the linear system satisfied by $\bs\phi$, namely
\begin{equation} \label{Num_sch-5}
     (D_1\,G_1+D_2\,G_2)\,\bs\phi^{k+1} = -\frac{2}{\delta t}\,\left(D_1\widetilde{\bs u}^{k+1}+D_2\widetilde{\bs v}^{k+1}\right).
\end{equation}
Once \eqref{Num_sch-5} is solved, the velocity is updated with \eqref{Num_sch-3}.
The linear system defined by \eqref{Num_sch-5} can be solved efficiently by
applying the same discrete transform used for the vertical velocity component.

The temperature equation \eqref{Boussinesq-1c-dt} is integrated over the volume cells $K_{i-\frac 12,j}$, leading to
\begin{align} \label{Num_sch-6}
 \bs\theta^{k+1} +\frac{\delta t}{2\sqrt{Ra\,Pr}}\, & A_2\,\bs\theta^{k+1}
                 = \bs\theta^k -\frac{\delta t}{2\sqrt{Ra\,Pr}}\,A_2\,\bs\theta^k  \nonumber \\
                & -\frac{\delta t}2\,\bs\psi_{\bs\gamma}\,\left(3\,N_3(\bs u^k,\bs v^k,\bs\theta^k)-N_3(\bs
                u^{k-1},\bs v^{k-1},\bs\theta^{k-1})\right)  \\
                & -\frac{\alpha\delta t}{2}\,\bs\psi_{\bs\gamma}\,(3\,\bs v^k-\bs v^{k-1}), \nonumber
\end{align}
where $\bs\psi_{\bs\gamma}=\{\psi_{\bs\gamma}^{ij}\}\in\mathbb{R}^{N(M-1)}$,
$\psi_{\bs\gamma}^{ij}=\psi_{\bs\gamma}(x_{i-\frac 12})$.
Boundary conditions for temperature are treated as for the vertical velocity component. 

\section{Numerical results}

The main purpose of this paper is to produce reference stationary solutions for heat island flows at Rayleigh numbers
$Ra\le 10^5$. In our study the stratification coefficient $\alpha$ is fixed as $\alpha=1$.
Dependency of solutions upon this parameter will be addressed in further works.

First, the accuracy of our code is evaluated by computing stationary solutions in a square 
differentially heated cavity. For this test case, benchmark solutions available in \cite{LeQuere} are used 
for comparison.
Then, stationary solutions for flows in a heat island are studied. The methodology used to produce accurate results
is detailed and solutions are described and analyzed.

\subsection{Validation of the code: the square differentially heated cavity test case}

In order to assess the validity of our code and to check the accuracy of the numerical scheme 
\eqref{Num_sch-1}-\eqref{Num_sch-6} we have performed numerical simulations of stationary 
solutions to the square differentially heated cavity for values of the Rayleigh number $Ra=10^6,\, 10^7$ and
$10^8$. 

Le Qu\'er\'e \cite{LeQuere} produced accurate benchmark solutions for this problem. In \cite{LeQuere}, Chebyshev 
polynomials were used for the spatial approximation and an influence matrix technique was applied in order to 
enforce the divergence free condition.
Note that this problem, described in \cite{DeVahl,LeQuere-DeRoquefort,LeQuere}, differs from the heat island problem 
by the computational domain and the boundary conditions. However, the discrete system
\eqref{Num_sch-1}-\eqref{Num_sch-6}, with $\gamma=0$ and $\alpha=0$, applies as well to this test case.

Stationary solutions were obtained on uniform grids, in both horizontal and vertical directions, with mesh sizes
decreasing from $1/32$ through $1/1024$. The choice of uniform grids is not optimal for this problem as boundary 
layers develop along the vertical heated walls: A large number of points is thus required in order to obtain accurate 
results. Such a choice is however convenient and allows us to easily check the code accuracy. 

The characteristic values suggested by De Vahl Davis in \cite{DeVahl} were computed and compared with those of 
the benchmark solutions \cite{LeQuere}. All these values are recovered and a second-order spatial convergence is
obtained (see Figure \ref{fig:conv-rate}).

\subsection{Stationary solutions of flows in a heat island}

Due to the presence of the vertical stratification, the thermal perturbations are convected in the horizontal 
direction far from the heated source line. 
As a consequence, very long domains have to be used in order to produce accurate results. 
Numerical simulations in small computational domains are contaminated by artificial boundary layers which develop
at the outflow boundaries $\vert x\vert=L/2$. If the domain length $L$ is not large enough, the temperature cannot 
smoothly relax towards the vertical stratified profile imposed on the boundaries. 

For fixed Rayleigh numbers and mesh sizes, numerical simulations in domains with increasing lengths are performed 
with the standard heat equation, that is \eqref{Boussinesq-1c-dt} with $\gamma=0$. 
The effects of $L$ and $H$ on the accuracy of the results are investigated. This approach while time consuming allows
to produce reference solutions. Numerical simulations with the truncated temperature equation ($\gamma=1$) are then 
performed for comparison. This study demonstrates the efficiency of the thermal sponge layers.
Finally, stationary solutions at Rayleigh numbers $Ra\le 10^5$ are computed on meshes with a finer resolution.
Characteristic values are reported and various profiles are reproduced and analyzed.

The stationary state of the numerical simulations was assumed to be reached when time variations of flow variables are
controlled as it follows
$$\max_{n\ge 0}\Big\{\frac{\vert\bs u^{n+1}-\bs u^n\vert_\infty}{\delta t},
                     \frac{\vert{\bs v^{n+1}-\bs v^n}\vert_\infty}{\delta t},
                     \frac{\vert{\bs \theta^{n+1}-\bs \theta^n}\vert_\infty}{\delta t}\Big\}\,\le\, \textrm{Tol}\ $$
where $\textrm{Tol}\in (10^{-10},10^{-8})$ is a given parameter.

\subsubsection{Numerical simulations in large elongated domains}

For fixed vertical resolutions, $h=1/16,\, 1/32$ and $1/64$, and Rayleigh numbers $Ra=10^3,\, 10^4$ and $10^5$,
numerical simulations were performed for increasing values of the sizes $L$ and $H$ of the computational domain. 
The discrete systems \eqref{Num_sch-1}-\eqref{Num_sch-6} with the parameter $\gamma$ set to zero was used which 
corresponds to the classical heat equation. The temperature fluctuation $\theta$ is the flow variable for which 
convergence with respect to the domain size is the slowest. Therefore, we choose to use as reference value to 
compare simulations in different domains the minimum value reached by $\theta$ inside the computational domain. 
For stationary solutions, the minimum is reached above the center of the heated element, namely on the axis $x=0$ 
and for $y>0$. Recalling that $\Omega=(-\frac L2,\frac L2)\times(0,H)$, we denote 
$\theta_\textrm{min}^{(L,H)}=\min_{(x,y)\in\Omega}{\theta(x,y)}$ 
and we define the reference values $(L_\textrm{ref},\, H_\textrm{ref})$ and $\theta^\textrm{ref}_\textrm{min}$ by
\begin{equation} \label{test_conv-1}
  \epsilon_{(L,H)}:=
  \frac{\vert\theta_\textrm{min}^{(L,H)}-\theta^\textrm{ref}_\textrm{min}\vert}{\vert\theta^\textrm{ref}_\textrm{min}\vert}\le
  0.1\, h^2 \textrm{ for } L\ge L_\textrm{ref},\, H\ge H_\textrm{ref}, 
\end{equation}
so that $\theta^\textrm{ref}_\textrm{min}$ is considered as a converged value for $\theta_\textrm{min}^{(L,H)}$. 
Values $(L_\textrm{ref},\, H_\textrm{ref})$ found for the considered Rayleigh numbers and mesh sizes are reported in
Table I. We observe that the domain length $L_\textrm{ref}$ is not sensitive to the Rayleigh number while the domain
height $H_\textrm{ref}$, for a fixed resolution, decreases when $Ra$ increases. 
The strength of the stratification grows with the Rayleigh number: The flow is pushed down to the ground.

Once these reference domain sizes have been found by repeated simulations in large domains, 
values of $(L,H)$ ensuring a $h^2$ approximation of $\theta_\textrm{min}^{(L,H)}$ were estimated. 
We introduce $(L_c,H_c)$ so that $\epsilon_{(L_c,H_c)}\le h^2$ and
\begin{equation} \label{test_conv-2}
  \epsilon_{(L,H)} \ge h^2 \textrm{ for } L<L_c,\, H<H_c. 
\end{equation}
Therefore, for fixed Rayleigh number and mesh size, $(L_c,H_c)$ are the minimum values of the domain sizes required to
compute a numerical solution accurate up to scheme accuracy. Estimates of these minimal values are reported in Table I. 
The tests \eqref{test_conv-1} and \eqref{test_conv-2} impose strong restrictions on the admissible domain
sizes. Indeed, errors are most often larger than $h^2$ even for second-order schemes (see Figure \ref{fig:conv-rate}
for example). The values $(L_c,H_c)$ estimated with \eqref{test_conv-2} are probably too restrictive. However, their
use ensures accurate results.

Figures \ref{fig:conv-history} represents the convergence history of $\epsilon_{(L,H)}$ with respect to the domain
length $L$ for $H=4, 6$ and $8$. The Rayleigh number is $Ra=10^4$ and the vertical resolution is $h=1/32$.
The convergence rate of $\theta_\textrm{min}^{(L,H)}$ towards the reference value $\theta^\textrm{ref}_\textrm{min}$ 
behaves like $1/L$ for small values of $L$ and like $L^{-3/2}$ for large values of $L$. The scheme accuracy is
reached for $H=6$ and $L$ of the order of $2\,000$ (see also Table I).
The same behavior with respect to $L$ was found for other vertical resolutions and Rayleigh numbers.

On Figure \ref{fig:conv-history-Ra}, the convergence history of $\epsilon_{(L,H)}$ for $H=H_c$
is represented for $Ra=10^3,\, 10^4$ and $10^5$. The convergence rate of $\epsilon_{(L,H)}$ is found to be
independent of the Rayleigh number: All curves have the same slope in logarithmic scales. Also, the value $h^2=1/32^2$
corresponding to the scheme accuracy is reached by $\epsilon_{(L,H_c)}$ for values of $L$ decreasing when $Ra$ is 
increased. 

For the vertical resolution $h=1/64$, the computation of the reference solution in the domain
$\Omega_\textrm{ref}=(0,L_\textrm{ref})\times(0,H_\textrm{ref})$ with $L_\textrm{ref}=10\,000$ 
was achieved with $28\, 000$ points in the horizontal direction. This guarantees that the mesh satisfies 
$\ell_i\approx h$ in the neighborhood of the heated element. 
With the use of the mapping function \eqref{map-mesh}, we have in that case: $\max{\ell_i}/\min{\ell_i}=74$.
Therefore $35.8$ (resp. $17.9$) millions of points were used to compute the reference solution at $Ra=10^3$ 
(resp. $10^5$) and $35\, 000$ (resp. $200\, 000$) time iterations were required to reach the stationary 
solution. This represents $600$ (resp. $1\,500$) monoprocessor computing hours on an IBM Power4 computer. 

Due to this need for a large amount of computational resources, such a study is not feasible on meshes with
smaller grid sizes in the vertical direction. By extrapolating the values $L_c$ and $H_c$ found for $h=1/16,1/32$ 
and $1/64$ (see Table I), we can roughly estimate these minimal values for a finer mesh, namely $h=1/128$. We obtain
$(L_c,H_c)=(12\,000,18)$ at $Ra=10^3$, $(L_c,H_c)=(10\,000,12)$ at $Ra=10^4$ and $(L_c,H_c)=(8\,000,10)$ at 
$Ra=10^5$. Therefore, a resolution of at least $115$ millions of points would be required in order to compute a 
reference solution on a grid with a vertical resolution $h=1/128$.
As it is shown in the next sections, the use of the truncated heat equation ($\gamma=1$) allows us to significantly relax 
these constraints on the computational parameters.

\subsubsection{Efficiency of the truncated temperature equation}

The effect of the filter function \eqref{psi-gamma} is to smoothly damp the convective terms in the heat equation 
in the vicinity of the domain boundaries $\vert x\vert=L/2$. By introducing such thermal sponge layers, we aim at 
improving the accuracy of the numerical simulations when the computational domains are not long enough to ensure 
a $O(h^2)$ approximation.

Numerical simulations of stationary solutions at $Ra=10^5$ have been performed for $H$ fixed to $H_c$, for the 
vertical resolutions listed in Table I and for increasing values of $L$. 
The errors $\epsilon(L,H_c)$ produced by the standard and the modified heat equation are used to compare the
efficiency of both models. 
The values $\sigma=0.85$ and $p=8.0$ have been retained for the filter function \eqref{psi-gamma}. 
These parameters were found to be efficient for the numerical simulations of stationary solutions.
A parametric study of the truncated temperature equation is beyond the scope of this paper. However, this 
question is important and will be addressed in future works on this problem. 

On Figure~\ref{fig:conv-trunc-history-1}, the errors $\epsilon(L,H_c)$ obtained for $h=1/32$ with the standard 
and the truncated heat equations are plotted. We note that:
\begin{description}
  \item[--] the truncated heat equation produces errors about $10$ times smaller than the standard equation even for small
            values of $L$;
  \item[--] both curves have the same decay rate and converge to the same asymptotic value;
  \item[--] the truncated heat equation produces values of $L_c$ which are approximately $3$ times smaller than the values
            listed in Table I and corresponding to the standard heat equation. 
            We recall that $L_c$, defined by \eqref{test_conv-2}, is the minimum value of the domain length required to
            ensure a $h^2$ approximation of the temperature fluctuation.
\end{description}
The values $L_c=120,\, 400$ and $900$ are obtained with the truncated equation for the respective resolutions $h=1/16,\,
1/32$ and $1/64$ while $L_c=480,\, 1\,280$ and $3\,200$ were necessary with the classical heat equation (see Table I). 

Also, by examining the time history of the discrete time variation 
$\frac{\vert{\bs \theta^{n+1}-\bs \theta^n}\vert_\infty}{\delta t}$, it appears that the convergence to the stationary
solution is achieved in less time iterations with the truncated equation than with the classical one. 
For example, in $\Omega=(-240,240)\times(0,4)$ and for $h=1/32$ the stationary solution at $Ra=10^5$ is reached
after $32\,000$ time iterations with the truncated heat equation while $48\,000$ time iterations are needed with
the classical one (see Figure \ref{fig:conv-theta}).

Therefore, for a given accuracy, stationary solutions can be computed with the truncated heat equation in 
significantly smaller computational domains than with the classical temperature equation and in less time iterations. 
This results in a use of less computational resources. As a consequence, this approach allows us to compute 
stationary solutions on meshes with a finer vertical resolution.

\subsubsection{Accurate stationary solutions} 
Direct numerical simulations at $Ra=10^3, 10^4$ and $10^5$ and with a vertical resolution $h=1/128$ have been
performed with the truncated temperature equation. The computational parameters are listed in Table II. 
The estimates derived in Section 5.2.1 are used to determine the computational domains. 
Also, in agreement with the previous section, domain lengths about $3$ times smaller than the estimated values are retained.

In order to characterize the stationary solutions, the maximum values of the velocity components $(u,v)$,
the temperature variation $\theta$, the vorticity $\omega=\partial v/\partial x-\partial u/\partial y$ and 
the streamfunction $\psi$ are reported in Table III. 
The locations in the computational domains where these extrema are reached are also collected: When only one of
the coordinates is listed, the other one is equal to $0$. 

Also reported in Table III is the Nusselt number which is defined by
$$ Nu\ = -\frac 1 L \int_{-L/2}^{L/2} \frac{\partial \theta}{\partial y}(x,y=0)\,dx. $$
As expected, this value, which measures the intensity of the heat transfer, increases with the Rayleigh number.

The velocity components and the temperature fluctuation decay rapidly with respect to the elevation $y$ as it is
shown on Figures \ref{fig:TV_profiles} and \ref{fig:U_profile} representing the vertical profiles at the center 
of the heated element $x=0$. The gravitational force and the vertical stratification limit the vertical propagation
of perturbations.
In the horizontal direction (see Figure \ref{fig:TV_x_profiles}b), the vertical velocity $v$
vanishes rapidly for $x$ outside of the heated region, that is for $\vert x\vert>0.5$. Therefore, 
vertical convection is essentially localized above the heated element: This behavior is independent 
of the Rayleigh number. However, its intensity increases with $Ra$. Indeed, the maximum value reached by
the vertical velocity component increases with $Ra$ (see Figure \ref{fig:TV_x_profiles}b and Table III). 
The temperature fluctuation and the horizontal velocity have a similar behavior in the horizontal direction for $\vert
x\vert\ge 2$ (see Figures \ref{fig:TV_x_profiles}b and \ref{fig:U_x_profile}): They decay slowly
to a small but nonzero value which is growing with $Ra$. Therefore, the convection outside the heated source line
is mainly horizontal. This illustrates the difficulty to approximate such flows in limited computational domains.

The profiles of the temperature fluctuation, plotted on Figures \ref{fig:TV_profiles}a and \ref{fig:TV_x_profiles}a,
show that the ascending propagation of the thermal perturbation is reduced when $Ra$ is increased. Simultaneously,
the profile of the velocity components exhibit largest extrema and steepest gradients for $y\le 1$ 
(see Figures \ref{fig:TV_profiles}b and \ref{fig:U_profile}). Also, elevations where the velocity
components are maximum decrease for growing $Ra$ (see Table III). Hence, when the Rayleigh number is increased, 
the flow is pushed down to the ground. 
At larger Rayleigh numbers, we expect that the competition between the natural convection, inducing an ascending
propagation, and the vertical stratification, limiting this effect, will induce a loss of symmetry of the solutions 
leading to unsteady flows.

To better illustrate the effect of the vertical stratification, isolines of the temperature fluctuation $\theta$ and
the vorticity $\omega$ are displayed on Figures \ref{fig:T_iso} and 
\ref{fig:Vort_iso} in a region surrounding the heat 
island perturbation, that is for $\vert x\vert\le 5$ and $y\le 3$. The thermal plume in form of a mushroom, 
typical in natural convection problems (see \cite{Xin-al} for instance), cannot develop in a stratified
medium. Instead, the main thermal structure is centered above the heated plate, symmetric with respect to 
the axis $x=0$ and very elongated in the horizontal direction. 
Above, a thermal sink characterized by negative temperature variation $\theta$ is observed. 
The intensity of this structure grows with $Ra$ while its vertical position decreases. 
The vorticity structures (see Figure \ref{fig:Vort_iso}) exhibit multi-cell symmetric patterns. 
They become thinner when $Ra$ increases and are clearly pushed down to the axis $y=0$.

\subsubsection{Computational efficiency}
A parallel version of the Fortran 90 code based on implicit communications (OpenMP) was used for the numerical simulations
presented in this paper: An efficiency of approximately $6.8$ is found on $8$ processors on a cluster of IBM Power 4 
computers. In order to perform the numerical simulations presented in Section 5.2.3, $6\,000$ monoprocessor hours were necessary. 
The CPU time per iteration and per node used by the code is $2.5\times 10^{-6}$ seconds on IBM Power 4 processors. 
Concerning the memory, $20$ real unknowns ($8$ bytes) have to be stored for each node of the mesh.

\section{Concluding remarks and perspectives}
In this paper, steady state solutions of a natural convection problem in an unbounded domain are investigated
by direct numerical simulations. For this problem, the flow is thermally stratified in the vertical 
direction and perturbed by a local heat island located on the ground. Due to the vertical stratification, 
the flow circulation is dominated by horizontal convection, so that perturbations are propagated in the 
horizontal direction far from the heated source.
Stationary solutions are first investigated by numerical simulations in \textit{very elongated} domains for
moderate vertical resolutions, that is $h=1/16,1/32$ and $1/64$. Repeated computations in increasing domains 
have been performed: the minimum length and height necessary to ensure a $O(h^2)$ accuracy have been estimated 
at $Ra=10^3,10^4$ and $10^5$. 
This approach, while time and memory consuming, provided reference simulations that have been used to
validate and compare results obtained with a truncated heat equation. We have shown that the use of a suitable thermal 
sponge layer placed at the vertical outflow allows to noticeably reduce the size of the 
computational domain. Therefore, numerical simulations on finer grids are made accessible. 
The stationary solutions at the aforementioned $Ra$ have been computed on grids with vertical resolution
$h=1/128$. Characteristic values of these steady states have been provided.

The thermal circulation induced by the heat island consists in symmetric multi-cell pattern centered above the heated
element. Flow structures are pushed down to the ground when the Rayleigh number is increased. Also, their intensity
grows with $Ra$. We therefore may expect that stability of steady states will be lost at larger $Ra$ leading 
to nonstationary solutions. The thermal sink found above the heat island should first oscillate with respect to 
the vertical axis $x=0$ in a periodic time regime. 
The numerical study of the development of instabilities and the detection of successive transitions from steady state
to turbulent flows is our main motivation. Contributions to this project will be presented in forthcoming papers.
Dependency of solutions upon the stratification coefficient is also an open question for this problem. Such study
will be addressed in future works.

\section*{Acknowledgments}
The numerical simulations presented in this paper were performed on the cluster of HP Proliant bi-processors of the 
Laboratoire de Math\'ematiques (Universit\'e Blaise Pascal and CNRS) and on the cluster of IBM Power 4 computers of the
Supercomputing Center IDRIS of CNRS (Orsay, France).

\newpage


\begin{figure}
\begin{pspicture}(14,5)

\psline[linewidth=1.25pt,arrowsize=6pt]{->}(2,1)(12,1)
\psline[linewidth=1.5pt]{|-|}(6.,1)(7.8,1)

\pscurve[linewidth=1.25pt,arrowsize=6pt]{->}(6.3,1.25)(6.1,2)(5.6,2.5)
\psline[linewidth=1.25pt,arrowsize=6pt]{->}(6.6,1.25)(6.5,2.)(6.3,2.5)
\psline[linewidth=1.25pt,arrowsize=6pt]{->}(6.9,1.25)(6.9,2.5)
\psline[linewidth=1.25pt,arrowsize=6pt]{->}(7.2,1.25)(7.3,2)(7.5,2.5)
\pscurve[linewidth=1.25pt,arrowsize=6pt]{->}(7.5,1.25)(7.7,2)(8.2,2.5)

\psline[linewidth=1.25pt,arrowsize=6pt]{->}(10,3.5)(10,2.5)
\uput{0.25}[0](10,3){$g$}

\uput{0.5}[-90](6.9,1){$0$}
\uput{0.5}[-90](6,1){$\frac\delta 2$}
\uput{0.5}[-90](5.75,1){$-$}
\uput{0.5}[-90](7.8,1){$\frac\delta 2$}

\end{pspicture}
\caption{Heat island perturbation.}
\label{fig:heat-island}
\end{figure}
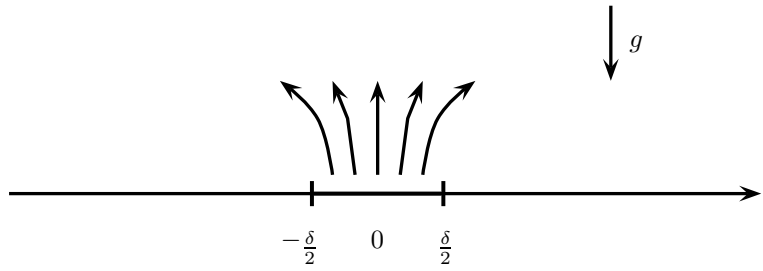

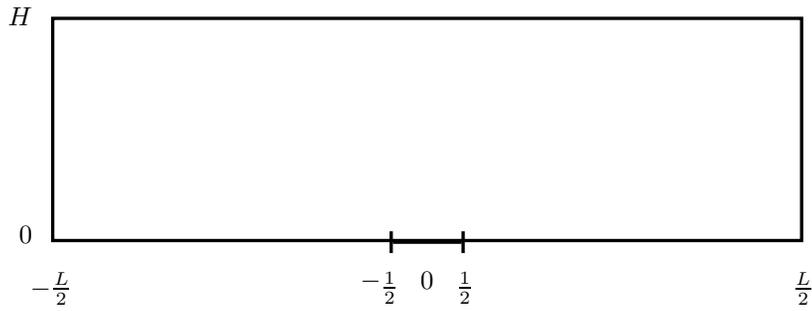
\begin{figure}
\begin{pspicture}(14,5)
\psframe[linewidth=1.25pt](2,1)(12,4)
\psline[linewidth=1.25pt]{|-|}(6.5,1)(7.5,1)
\uput{0.25}[-180](2,1.1){$0$}
\uput{0.25}[-180](2,4){$H$}
\uput{0.4}[-90](2,1){$-\frac L 2$}
\uput{0.4}[-90](12,1){$\frac L 2$}
\uput{0.4}[-90](7,1){$0$}
\uput{0.4}[-90](6.5,1){$\frac 12$}
\uput{0.4}[-90](6.25,1){$-$}
\uput{0.4}[-90](7.5,1){$\frac 12$}
\end{pspicture}
\caption{Computational domain $\Omega=\bigl(-\frac L2,\frac L2\bigr)\times(0,H)$.}
\label{fig:comp-domain}
\end{figure}

\begin{figure}
\begin{pspicture}(14,5)

\psframe[linewidth=1.5pt](5.5,1)(8.5,3)
\uput{0.5}[-90](5.5,1){$x_{i-1}$}
\uput{0.5}[-90](7,1){$x_{i-1/2}$}
\uput{0.5}[-90](8.5,1){$x_i$}
\uput{0.5}[-90](10,1){$x_{i+1/2}$}
\uput{0.5}[-180](5.5,1){$y_{j-1}$}
\uput{0.5}[-180](5.5,2){$y_{j-1/2}$}
\uput{0.5}[-180](5.5,3){$y_j$}
\uput{0.5}[-180](5.5,4){$y_{j+1/2}$}
\qdisk(7,2){3pt}
\uput[-135](7,2){$P_{ij}$}
\pspolygon[fillstyle=solid,fillcolor=black](8.4,1.9)(8.4,2.1)(8.6,2)
\uput[0](8.6,2){$u_{ij}$}
\pspolygon[fillstyle=solid,fillcolor=black](6.9,2.9)(7.1,2.9)(7,3.1)
\uput[45](7,3.1){$v_{ij}$}
\uput[135](7,3.1){$\theta_{ij}$}
\psframe[linecolor=gray,linewidth=1pt,linestyle=dashed](5.5,2)(8.5,4)
\psframe[linecolor=darkgray,linewidth=1pt,linestyle=dotted](7,1)(10,3)
\end{pspicture}
\caption{Cells $K_{i-\frac 12,j-\frac 12}$ (solid), $K_{i-\frac 12,j}$ (dashed) and $K_{i,j-\frac 12}$ (dotted) and
their corresponding discrete values.}
\label{fig-cell}
\end{figure}
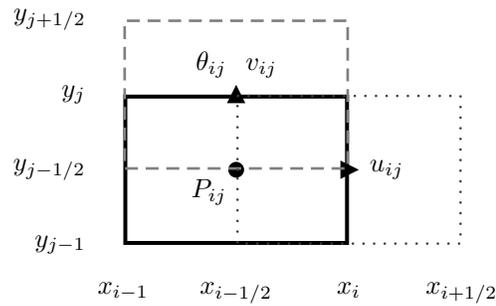


\begin{figure}
\begin{center}
\vlegend{1.}{Relative error}{.1} \hspace{.5cm}
\includegraphics[width=4.in]{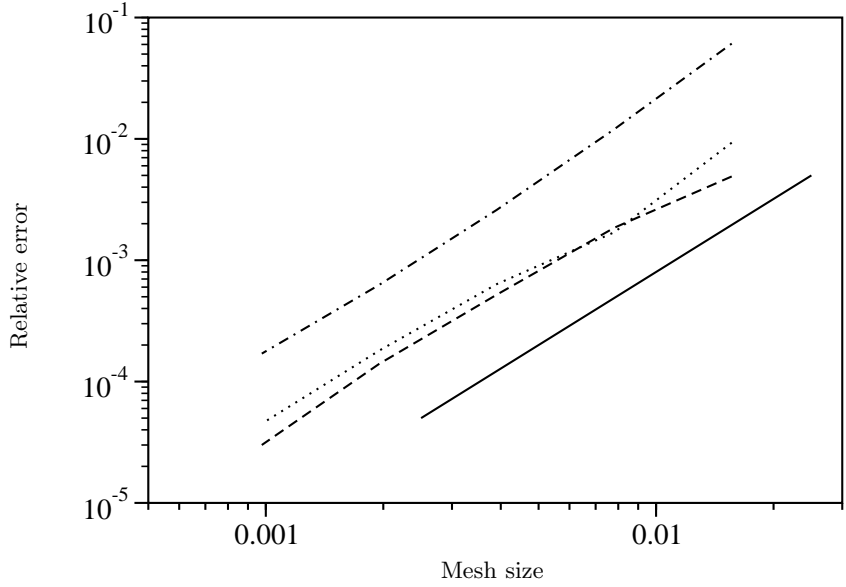} \\
\hspace{1.5cm}\small{Mesh size}
\end{center}
\caption{Spatial accuracy for the square differentially heated cavity test case at $Ra=10^6$. 
         Solid line: slope 2 ; dashed line: horizontal velocity ; dotted line: vertical velocity ; 
         dashed-dotted line: Nusselt number.}
\label{fig:conv-rate}
\end{figure}

\begin{table}
\caption{Estimate of the size of the computational domain for various mesh size $h$ and Rayleigh number $Ra$.}
\begin{center}
\begin{small}
\begin{tabular}{ccccc}
\hline
     $Ra$      &  $N_y/H$   &  $(L_c,H_c)$       & $(L_{\textrm{ref}},H_{\textrm{ref}})$   \\
\hline
  $10^3$        &    $16$    &   $(960,8)$        &   $(3\,200,12)$         \\
                &    $32$    &   $(2\,560,10)$    &   $(6\,200,16)$         \\
                &    $64$    &   $(5\,200,14)$    &   $(10\,000,20)$         \\[0.05in]
  $10^4$        &    $16$    &   $(640,5)$        &   $(3\,200,8)$          \\
                &    $32$    &   $(1\,920,6)$     &   $(6\,200,10)$         \\
                &    $64$    &   $(4\,200,8)$     &   $(10\,000,12)$        \\[0.05in]
  $10^5$        &    $16$    &   $(480,3)$        &   $(3\,200,6)$          \\
                &    $32$    &   $(1\,280,4)$     &   $(6\,200,8)$          \\
                &    $64$    &   $(3\,200,6)$     &   $(10\,000,10)$        \\ 
\hline
\end{tabular}
\end{small}
\end{center}
\end{table}

\begin{figure}
\begin{center}
\vlegend{.9}{Relative error}{.1} \hspace{.15cm}
\includegraphics[width=3.9in]{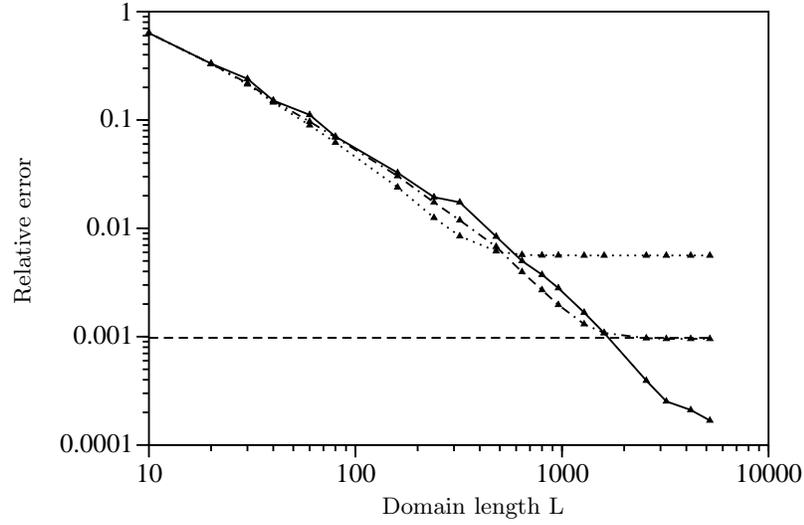} \\
\vspace{0.cm}\hspace{1.5cm}\small{Domain length L}
\end{center}
\caption{Convergence of $\epsilon_{(L,H)}$ at $Ra=10^4$ with respect to the domain length for $H=4$ 
         (dotted line), $H=6$ (dotted-dashed line) and $H=8$ (solid line).
         The dashed line corresponds to the expected accuracy $h^2=1/32^2$. }
\label{fig:conv-history}
\end{figure}

\begin{figure}
\begin{center}
\vlegend{.9}{Relative error}{.1} \hspace{.15cm}
\includegraphics[width=3.9in]{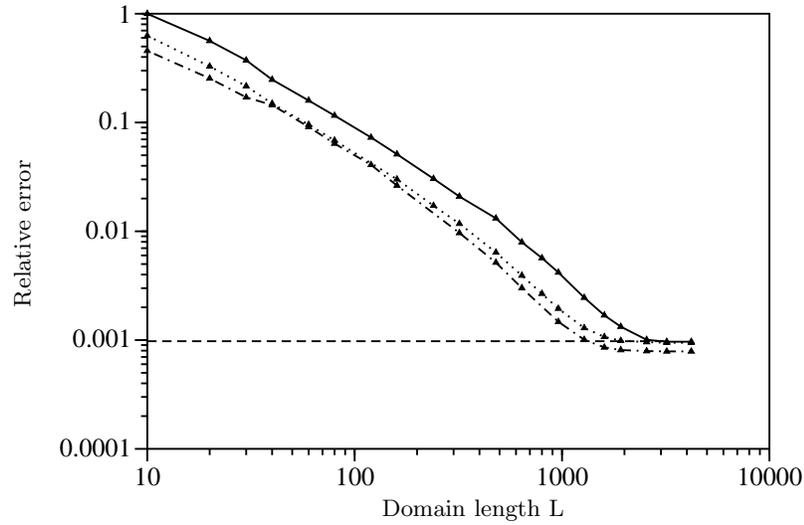} \\
\vspace{0.cm}\hspace{1.5cm}\small{Domain length L}
\end{center}
\caption{Convergence of $\epsilon_{(L,H_c)}$ with respect to the domain length for $Ra=10^3$ 
         (solid line), $Ra=10^4$ (dotted line) and $Ra=10^5$ (dashed-dotted line).
         The dashed line corresponds to the expected accuracy $h^2=1/32^2$. }
\label{fig:conv-history-Ra}
\end{figure}

\begin{figure}
\begin{center}
\vlegend{.9}{Relative error}{.1} \hspace{.15cm}
\includegraphics[width=3.8in]{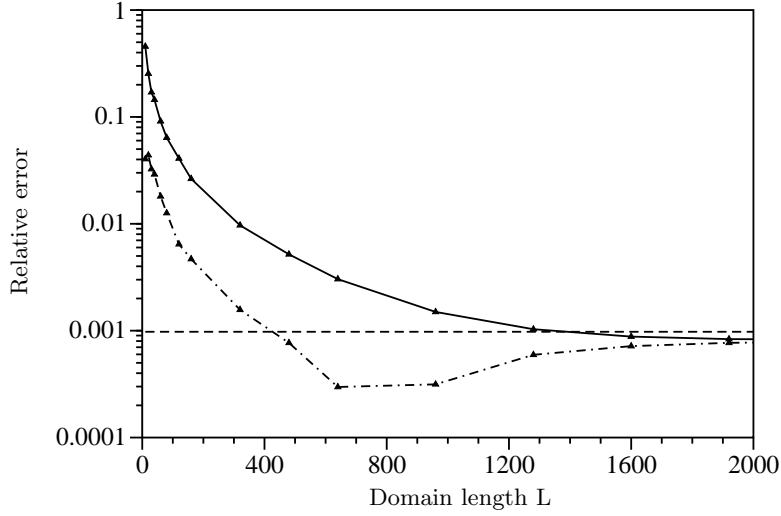} \\
\vspace{0.cm}\hspace{1.5cm}\small{Domain length L}
\end{center}
\caption{Convergence of $\epsilon_{(L,H=4)}$ at $Ra=10^5$ with respect to the domain length obtained 
         with the standard heat equation (solid line) and the truncated version (dashed-dotted line).
         The dashed line corresponds to the expected accuracy $h^2=1/32^2$. }
\label{fig:conv-trunc-history-1}
\end{figure}

\begin{figure}
\begin{center}
\vlegend{.9}{Time variations}{.1} \hspace{.15cm}
\includegraphics[width=3.8in]{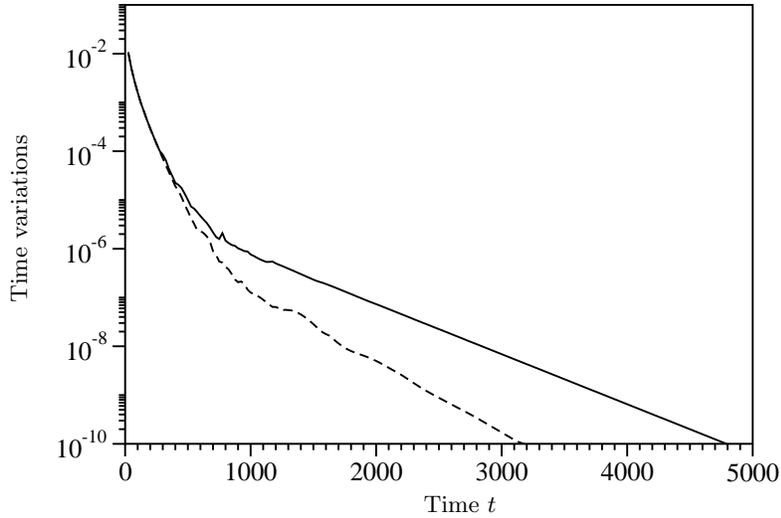} \\
\vspace{0.cm}\hspace{1.5cm}\small{Time $t$}
\end{center}
\caption{Convergence history of $\max_{i,j}\frac{\vert\theta_{i,j}^{n+1}-\theta_{i,j}^n\vert}{\delta t}$ for the standard 
         (solid line) and the truncated (dotted line) heat equation at $Ra=10^5$. The vertical resolution is
         $h=1/32$, the time step is $\delta t=0.1$ and the computational domain is $\Omega=(-240,240)\times(0,4)$.}
\label{fig:conv-theta}
\end{figure}

\begin{table}
\caption{Computational parameters for numerical simulations on meshes with $h=1/128$.}
\begin{center}
\begin{small}
\begin{tabular}{cccc}
\hline
 $Ra$                     &  $10^3$         &  $10^4$           & $10^5$           \\
\hline
  $(L,H)$                  & $(4\,400,18)$       &   $(3\,600,12)$        & $(2\,560,10)$         \\
  $(N,M)$                  & $(19\,000,2\,304)$  &   $(16\,500,1\,536)$   & $(12\,000,1\,280)$    \\
  $\delta t$               & $0.1$               &   $0.05$               & $0.025$               \\
  $T_\textrm{stat}$        & $4\,893$            &   $6\,569$             & $10\,902$             \\
\hline
\end{tabular}
\end{small}
\end{center}
\end{table}

\begin{table}
\caption{Characteristic values of the stationary solutions.}
\begin{center}
\begin{small}
\begin{tabular}{cccc}
\hline
 $Ra$                    &      $10^3$              &       $10^4$             &      $10^5$                 \\
\hline
  $\theta_{\min}$         &  $-0.024823$             &  $-0.071289$             &    $-0.166316$              \\
   $y$                    &   $1.36191$              &   $0.94640$              &     $0.84320$               \\[0.05in]
$u_{\max}$                &   $0.118887$             &   $0.174844$             &     $0.179054$              \\
   $(x,y)$                &  $(-0.52314,0.17705)$    & $(-0.36781,0.13255)$     & $(-0.30643,0.09322)$        \\[0.05in]
$v_{\max}$                &   $0.125594$             &   $0.228250$             &     $0.322483$              \\
   $y$                    &   $0.55124$              &   $0.44684$              &     $0.42552$               \\[0.05in]
$v_{\min}$                &  $-0.030470$             &  $-0.039291$             &    $-0.079265$              \\
   $(x,y)$                &  $(-0.85107,0.41913)$    & $(-0.60361,0.32122)$     &  $(-0.39114,0.64777)$        \\[0.05in]
$\omega_{\max}$           &   $2.06900$              &    $3.951325$            &     $5.921375$              \\
   $x$                    &   $0.49642$              &    $0.48402$             &     $0.47813$               \\[0.05in]
$\psi_{\max}$             &   $0.042954$             &    $0.048265$            &     $0.040272$              \\
   $(x,y)$                &  $(-0.60014,0.59882)$    & $(-0.38883,0.45575)$     & $(-0.25610,0.45867)$        \\[0.05in]
$Nu$                      &   $0.148605$             &    $0.295132$            &     $0.643594$              \\
\hline
\end{tabular}
\end{small}
\end{center}
\end{table}

\begin{figure}
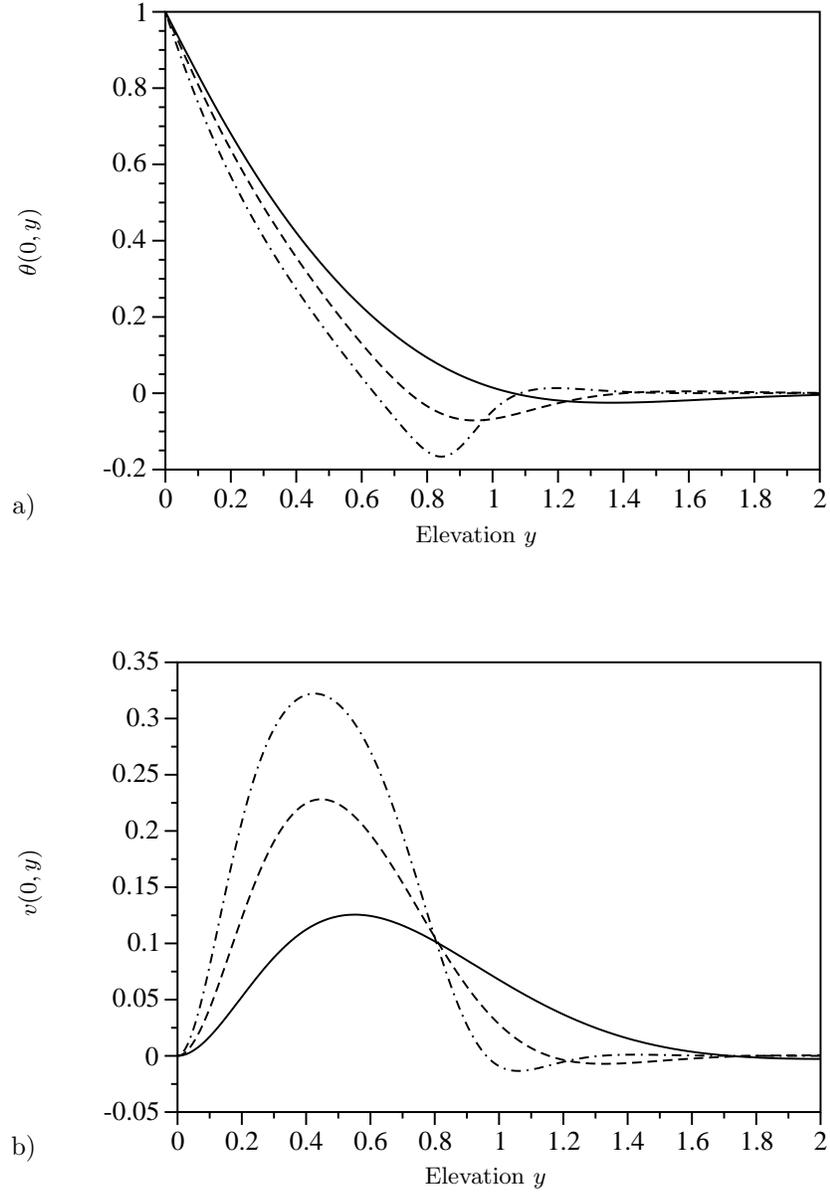

\begin{center}
a)\vlegend{1.2}{$\theta(0,y)$}{-.1} \hspace{.15cm}
\includegraphics[width=3.8in]{pft.eps} \\
\vspace{0.cm}\hspace{1.5cm}\small{Elevation $y$}
\end{center}
\vspace{0.25in}
\begin{center}
b)\vlegend{1.2}{$v(0,y)$}{-.1} \hspace{.15cm}
\includegraphics[width=3.8in]{pfv.eps} \\
\vspace{0.cm}\hspace{1.75cm}\small{Elevation $y$}
\end{center}
\caption{Profiles of the temperature variation $\theta$ (a) and of the vertical velocity $v$ 
         (b) at the center of the heated element, \textit{i.e.} at $x=0$, for $Ra=10^3$ (solid line), 
         $Ra=10^4$ (dashed line) and $Ra=10^5$ (dashed-dotted line). 
         The vertical resolution $h=1/128$ and the truncated heat equation were used.}
\label{fig:TV_profiles}
\end{figure}

\begin{figure}
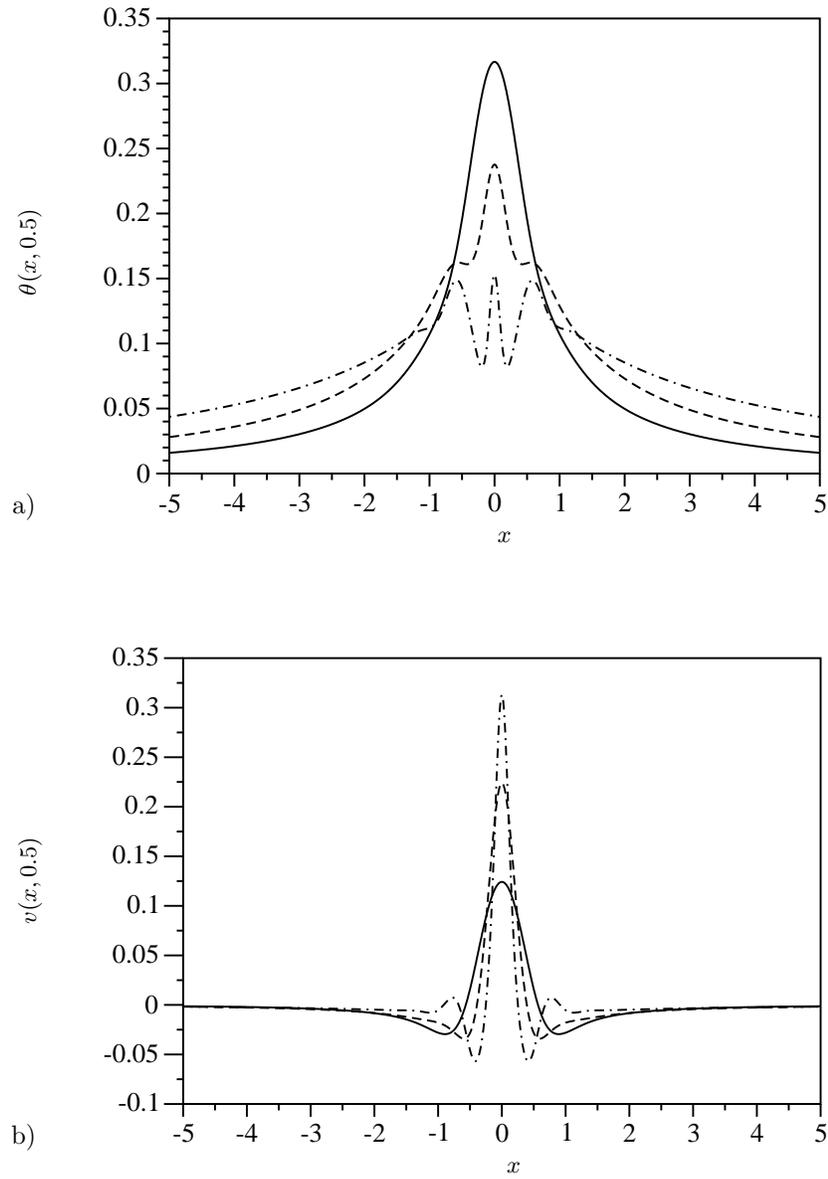

\begin{center}
a)\vlegend{1.1}{$\theta(x,0.5)$}{-.1} \hspace{.15cm}
\includegraphics[width=3.8in]{pfx_t.eps} \\
\vspace{0.cm}\hspace{2.25cm}\small{$x$}
\end{center}
\vspace{0.25in}
\begin{center}
b)\vlegend{1.1}{$v(x,0.5)$}{-.1} \hspace{.15cm}
\includegraphics[width=3.8in]{pfx_v.eps} \\
\vspace{0.cm}\hspace{2.5cm}\small{$x$}
\end{center}
\caption{Profiles of the temperature variation $\theta$ (a) and of the vertical velocity $v$ 
         (b) at the elevation $y=0.5$, for $Ra=10^3$ (solid line), $Ra=10^4$ (dashed line) and 
         $Ra=10^5$ (dashed-dotted line). 
         The vertical resolution $h=1/128$ and the truncated heat equation were used.}
\label{fig:TV_x_profiles}
\end{figure}

\begin{figure}
\begin{center}
\vlegend{1.1}{$u(0.25,y)$}{-.1} \hspace{.15cm}
\includegraphics[width=3.7in]{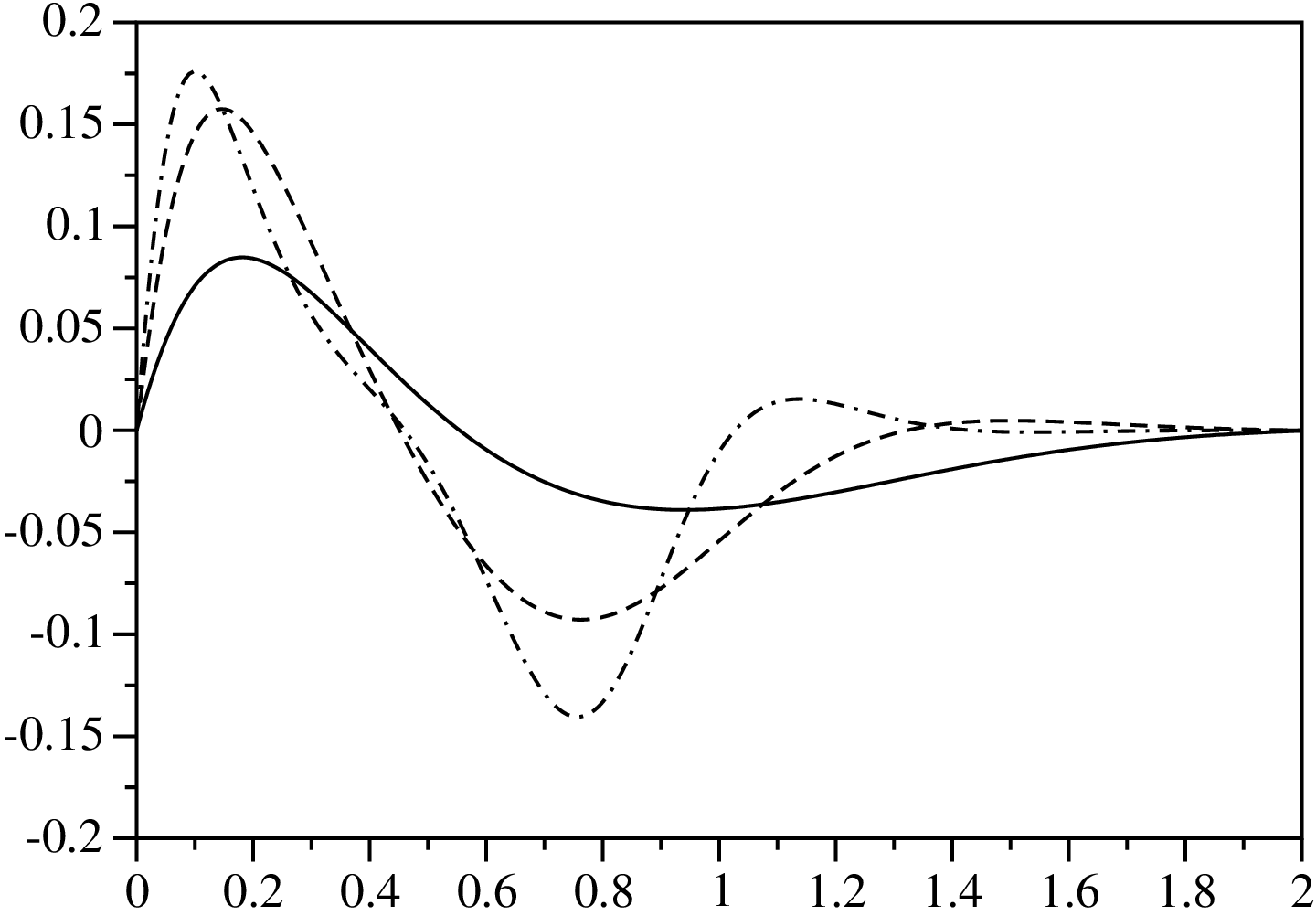} \\
\vspace{0.cm}\hspace{1.5cm}\small{Elevation $y$}
\end{center}
\caption{Profile of the horizontal velocity $u$ at $x=0.25$ for $Ra=10^3$ (solid line), 
         $Ra=10^4$ (dotted line) and $Ra=10^5$ (dashed line). 
         The vertical resolution $h=1/128$ and the truncated heat equation were used.}
\label{fig:U_profile}
\end{figure}

\begin{figure}
\begin{center}
\vlegend{1.1}{$u(x,0.1)$}{-.1} \hspace{.15cm}
\includegraphics[width=3.7in]{pfx_u.eps} \\
\vspace{0.cm}\hspace{2.cm}\small{$x$}
\end{center}
\caption{Profile of the horizontal velocity $u$ at the elevation $y=0.1$ for $Ra=10^3$ (solid line), 
         $Ra=10^4$ (dotted line) and $Ra=10^5$ (dashed line). 
         The vertical resolution $h=1/128$ and the truncated heat equation were used.}
\label{fig:U_x_profile}
\end{figure}

\begin{figure}
\vspace{0.3cm}
\begin{center}
\includegraphics[width=5.25in]{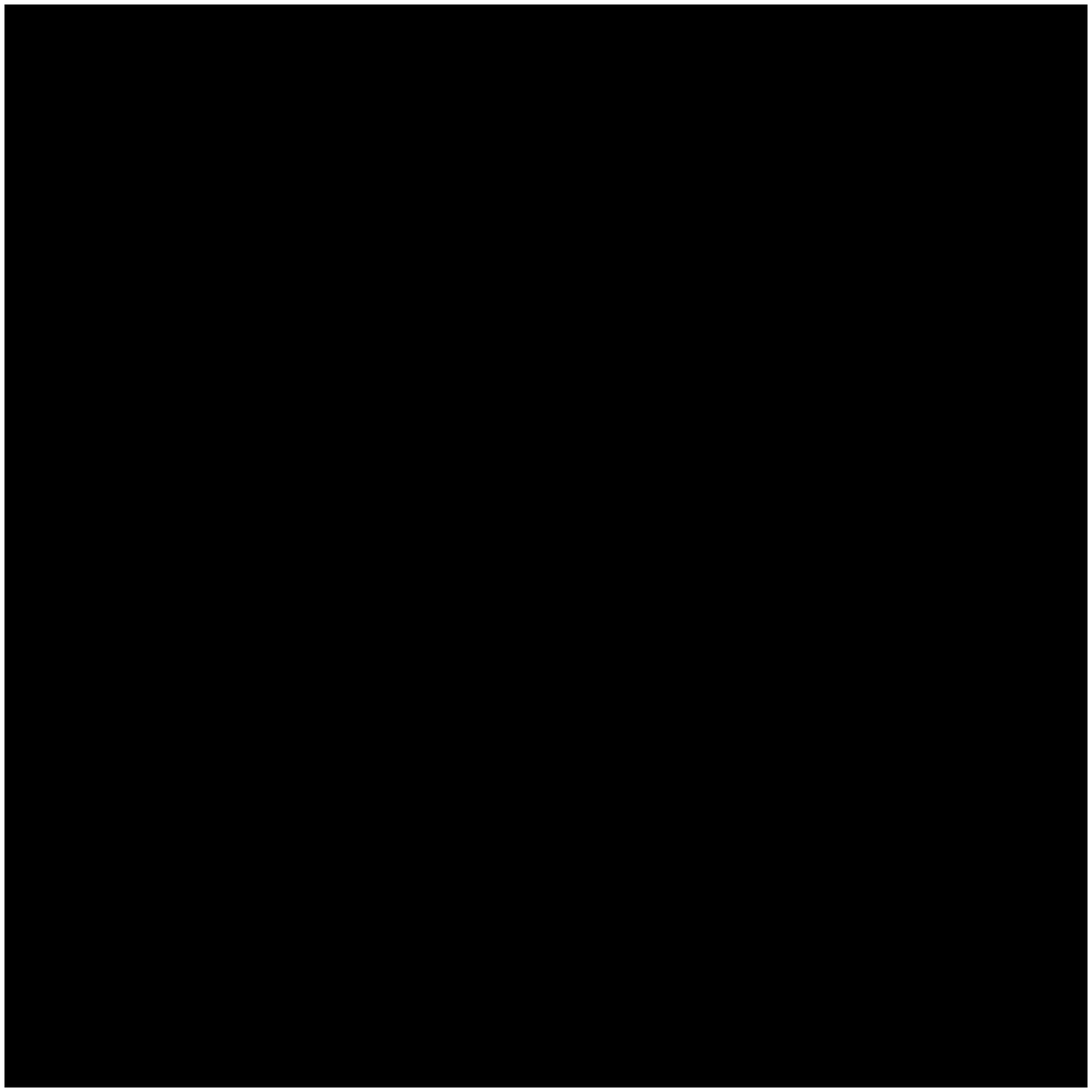} \\
\end{center}
\vlegend{0.7}{Elevation $y$}{-1.5}  \\[-0.5cm]
$Ra=10^3$\hspace{5.75cm}\small{$x$} \\[0.2cm]
Dashed lines: $-0.024,-0.02,-0.015,-0.01,-0.008,-0.007,-0.006,-0.005,-0.004,-0.003,-0.002,$ $-0.001.$
\begin{center}
\includegraphics[width=5.25in]{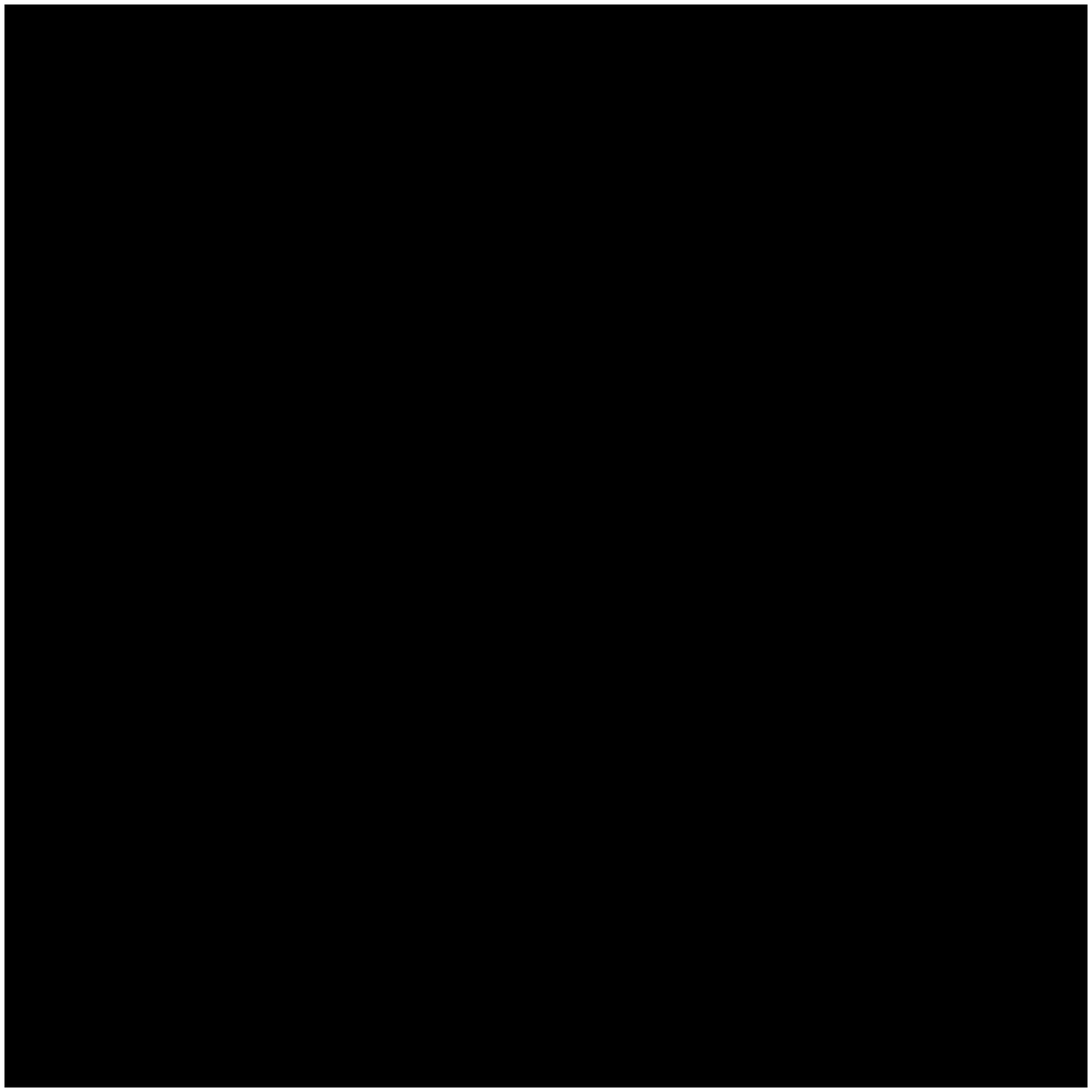} \\
\end{center}
\vlegend{0.7}{Elevation $y$}{-1.5}  \\[-0.5cm]
$Ra=10^4$\hspace{5.8cm}\small{$x$} \\[0.2cm]
Dashed lines: $-0.07,-0.06,-0.04,-0.02,-0.011,-0.01,-0.008,-0.006,$ $-0.004$, $-0.002,-0.001.$
\begin{center}
\includegraphics[width=5.25in]{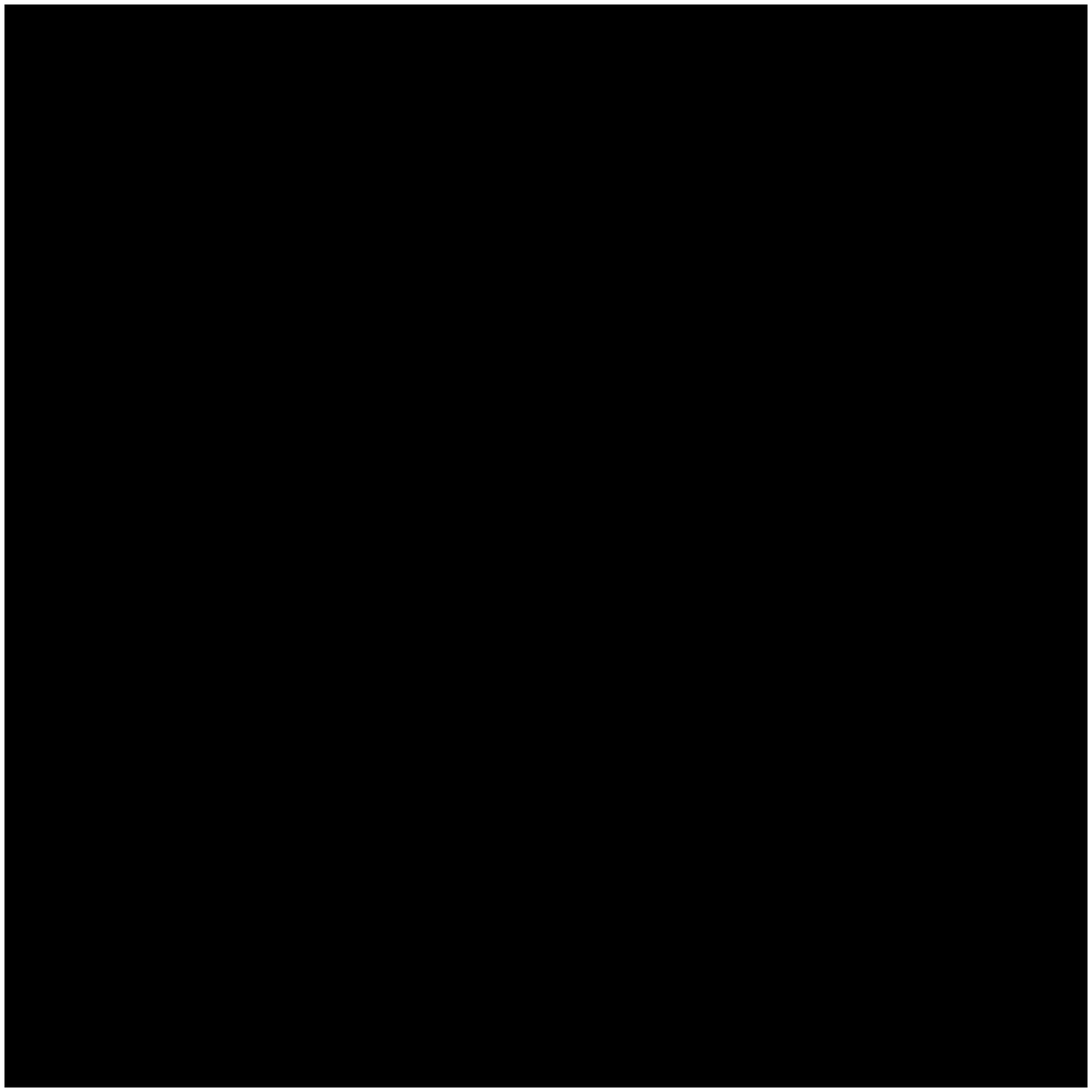} \\
\end{center}
\vlegend{0.7}{Elevation $y$}{-1.5}  \\[-0.5cm]
$Ra=10^5$\hspace{5.8cm}\small{$x$} \\[0.2cm]
Dashed lines: $-0.16,-0.13,-0.1,-0.06,-0.02,-0.016,-0.014,-0.012,-0.01,-0.008,-0.006,-0.004,$ $-0.002.$
\caption{Isolines of the temperature variation for increasing Rayleigh numbers. 
         Solid lines: $0.01$, $0.02$, $0.03$, $0.04$, $0.05$, $0.06$, $0.07$, $0.08$,
         $0.09$, $0.1$, $0.12$, $0.14$, $0.16$, $0.2$, $0.3$, $0.5$, $0.7$.}
\label{fig:T_iso}
\end{figure}

\begin{figure}
\begin{center}
\includegraphics[width=5.25in]{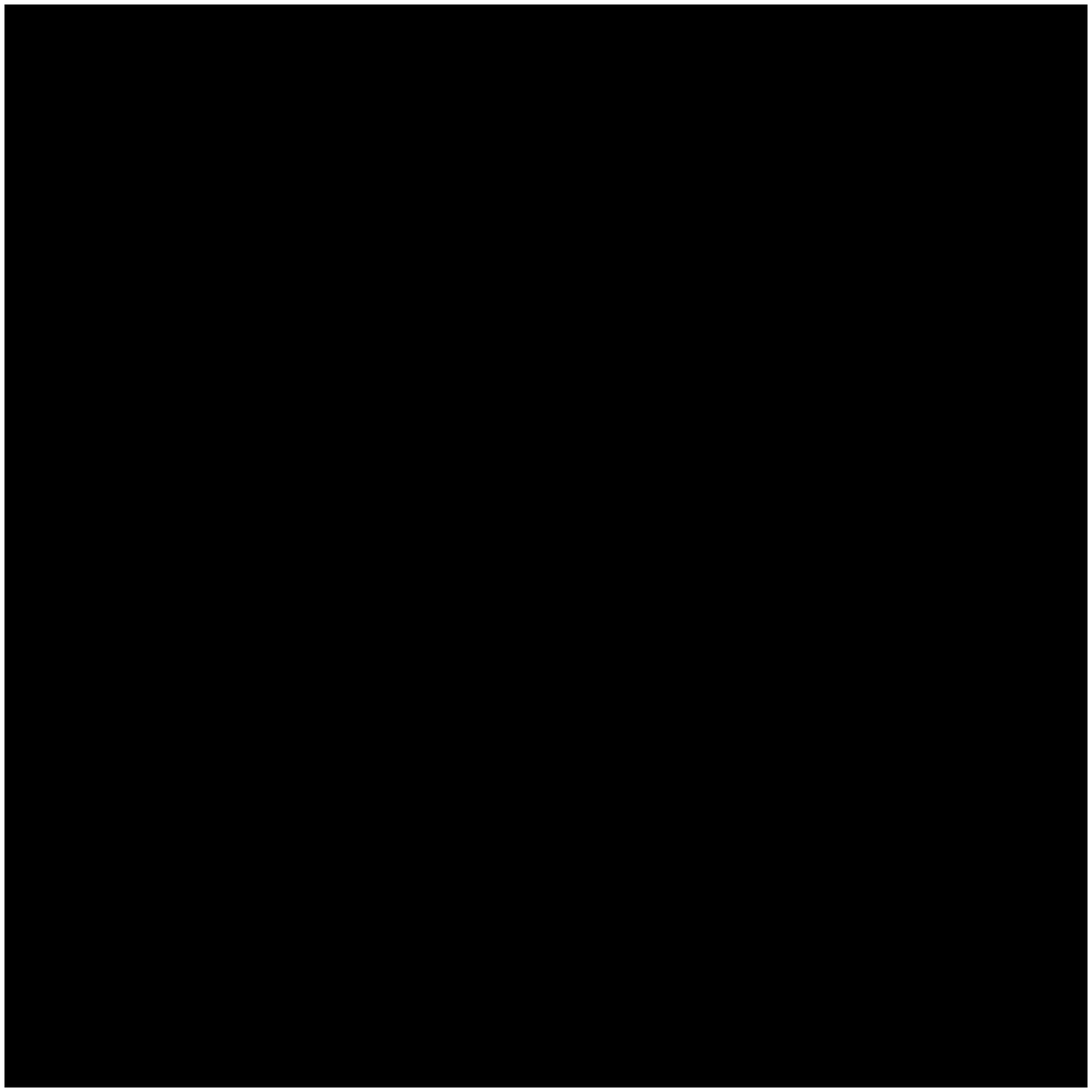} \\
\end{center}
\vlegend{0.7}{Elevation $y$}{0.5} \\[-0.5cm]
$Ra=10^3$ \hspace{5.6cm}\small{$x$} \\[0.2cm]
Solid lines: $0.02,0.03,0.04,0.05,0.06,0.07,0.078,0.1,0.15,0.2,0.3,0.4,0.5,0.6,1.0.$ \\
Dashed lines are opposite values
\begin{center}
\includegraphics[width=5.25in]{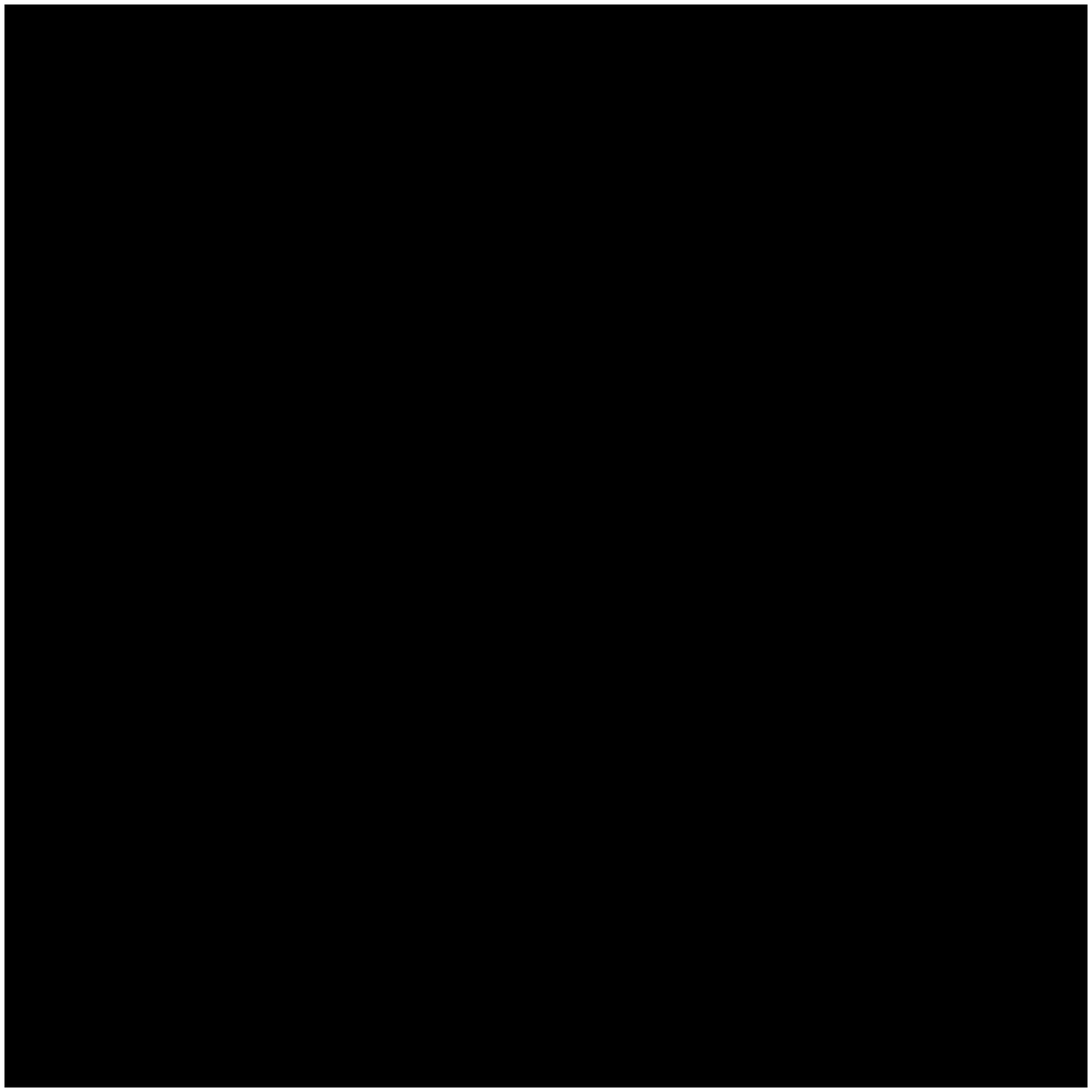} \\
\end{center}
\vlegend{0.7}{Elevation $y$}{0.5} \\[-0.5cm]
$Ra=10^4$ \hspace{5.75cm}\small{$x$} \\[0.2cm]
Solid lines: $0.02,0.04,0.06,0.08,0.1,0.15,0.2,0.25,0.3,0.4,0.6,1.0$ \\
Dashed lines are opposite values.
\begin{center}
\includegraphics[width=5.25in]{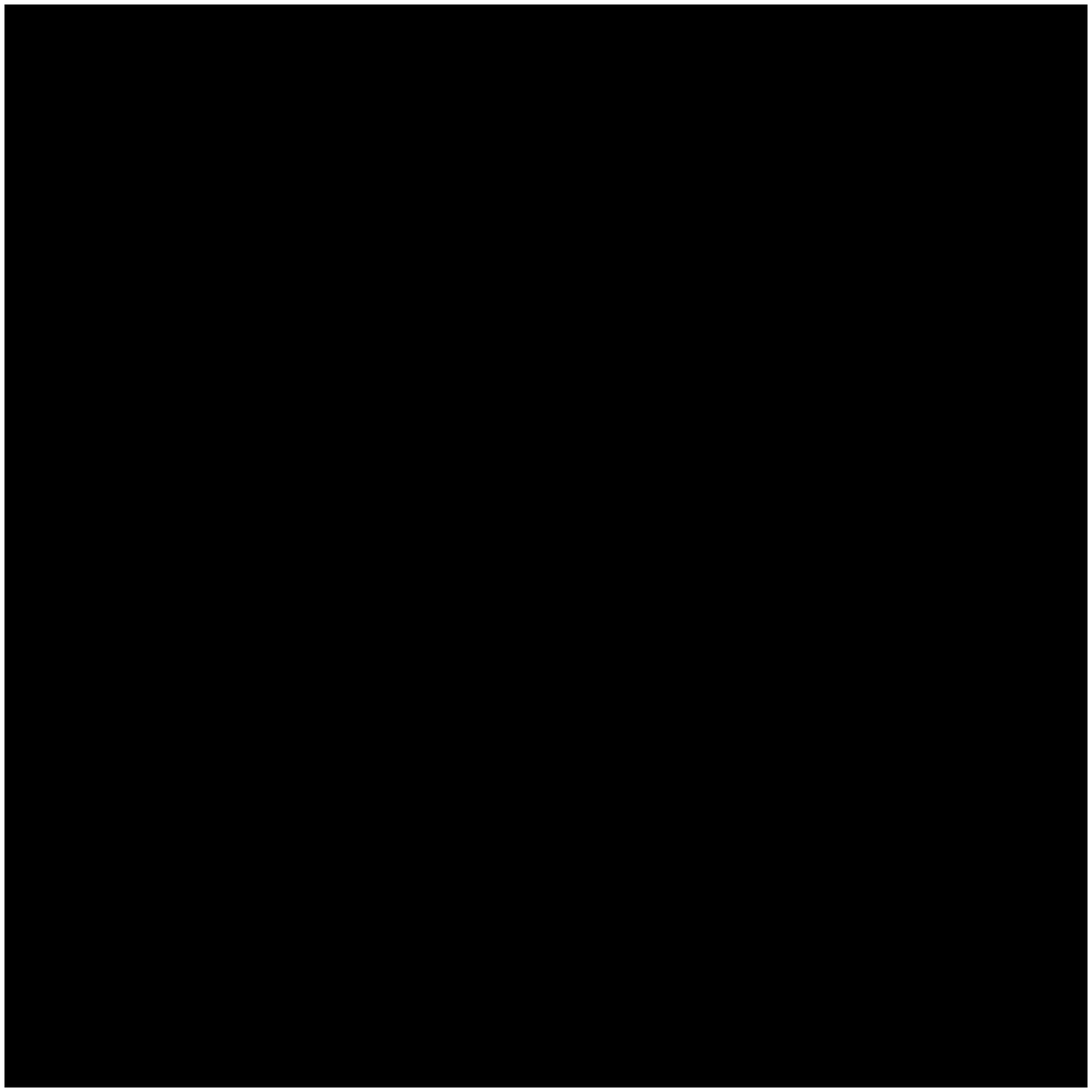} \\
\end{center}
\vlegend{0.7}{Elevation $y$}{0.5} \\[-0.5cm]
$Ra=10^5$ \hspace{5.75cm}\small{$x$} \\[0.2cm]
Solid lines: $0.01,0.02,0.04,0.06,0.07,0.1,0.15,0.2,0.3,0.4,0.6,1,1.5,2.0.$ \\
Dashed lines are opposite values.
\caption{Isolines of the vorticity for increasing Rayleigh numbers.}
\label{fig:Vort_iso}
\end{figure}

\end{document}